\input amstex
\documentstyle{amsppt}
\document

\magnification 1100

\def\gen{{\frak{g}}}

\def\a{{\alpha}}
\def\g{{\gamma}}
\def\o{{\omega}}

\def\l{{\lambda}}
\def\b{{\beta}}
\def\eps{{\varepsilon}}

\def\1b{{\bold 1}}

\def\cb{{\bold c}}

\def\mb{{\bold m}}

\def\Ab{{\bold A}}

\def\Bb{{\bold B}}

\def\Db{{\bold D}}

\def\Ib{{\bold I}}

\def\Kb{{\bold K}}
\def\Hb{{\bold H}}

\def\Rb{{\bold R}}
\def\Sb{{\bold S}}
\def\Tb{{\bold T}}

\def\d{{\roman d}}
\def\e{{\roman e}}

\def\k{{\roman k}}

\def\M{{\roman M}}

\def\A{{\roman A}}

\def\F{{\roman F}}

\def\R{{\roman R}}

\def\K{{\roman K}}

\def\Hbu{{\underline{\Hb}}}

\def\Scbu{{\underline{\Sb\cb}}}
\def\Scu{{\underline{\Sc}}}
\def\Ccu{{\underline{\Cc}}}
\def\Tbu{{\underline{\Tb}}}

\def\Du{{\underline{D}}}

\def\Pu{{\underline{P}}}

\def\Vu{{\underline{V}}}

\def\Dcu{{\underline{\Dc}}}
\def\Ocu{{\underline{\Oc}}}

\def\ad{{\roman{ad}}}

\def\Spec{{\roman{Spec}}}
\def\Aut{\text{Aut}}

\def\Hom{{\roman{Hom}}}

\def\dim{{\roman{dim}}}

\def\End{{\roman{End}}}
\def\ch{{\roman{ch}}}
\def\Id{{\roman{Id}}}

\def\Oplus{\ts\bigoplus}
\def\Sum{\ts\sum}
\def\Prod{\ts\prod}

\def\Ker{\text{Ker}\,}

\def\GL{{\roman{GL}}}
\def\SL{{\roman{SL}}}

\def\re{{\roman{re}}}
\def\reg{{\circ}}

\def\CC{{\Bbb C}}

\def\QQ{{\Bbb Q}}
\def\RR{{\Bbb R}}

\def\ZZ{{\Bbb Z}}

\def\Ac{{\Cal A}}
\def\Bc{{\Cal B}}
\def\Cc{{\Cal C}}
\def\Dc{{\Cal D}}

\def\Hc{{\Cal H}}
\def\Ic{{\Cal I}}

\def\Lc{{\Cal L}}
\def\Mc{{\Cal M}}
\def\Nc{{\Cal N}}
\def\Oc{{\Cal O}}
\def\Pc{{\Cal P}}

\def\Sc{{\Cal S}}

\def\Modc{{{{\Mc od}}}}
\def\mod{{{{mod}}}}

\def\mof{{{{mof}}}}

\def\and{{\quad\text{and}\quad}}

\def\ts{\textstyle}

\def\qed{\hfill $\sqcap \hskip-6.5pt \sqcup$}        
\overfullrule=0pt                                    

\def\pro{{\lim\limits_{\lla}}}
\def\ind{{\lim\limits_{\lra}}}

\def\la{{\langle}}
\def\ra{{\rangle}}
\def\lla{{\longleftarrow}}
\def\lra{{{\longrightarrow}}}

\newdimen\Squaresize\Squaresize=14pt
\newdimen\Thickness\Thickness=0.5pt
\def\Square#1{\hbox{\vrule width\Thickness
	      \alphaox to \Squaresize{\hrule height \Thickness\vss
	      \hbox to \Squaresize{\hss#1\hss}
	      \vss\hrule height\Thickness}
	      \unskip\vrule width \Thickness}
	      \kern-\Thickness}
\def\Vsquare#1{\alphaox{\Square{$#1$}}\kern-\Thickness}

\topmatter
\title From double affine Hecke algebras to quantized affine Schur algebras
\endtitle
\rightheadtext{}
\abstract We prove an equivalence between some category of modules of
the double affine Hecke algebra of type $A$
and of the quantized affine Schur algebra.
\endabstract
\author M. Varagnolo, E. Vasserot\endauthor
\address D\'epartement de Math\'ematiques,
Universit\'e de Cergy-Pontoise, 2 av. A. Chauvin,
BP 222, 95302 Cergy-Pontoise Cedex, France\endaddress
\email eric.vasserot\@math.u-cergy.fr\endemail
\address D\'epartement de Math\'ematiques,
Universit\'e de Cergy-Pontoise, 2 av. A. Chauvin,
BP 222, 95302 Cergy-Pontoise Cedex, France\endaddress
\email michela.varagnolo\@math.u-cergy.fr\endemail
\thanks
2000{\it Mathematics Subject Classification.}
Primary 17B37; Secondary 17B67, 14M15, 16E20.
\endthanks
\endtopmatter
\document

\head Introduction\endhead

Let $\F$ be a local non Archimedian field of residual characteristic $p$,
$q$ the order of the residual field.
Let $\k$ be an algebraically closed field of characteristic $\ell$.
Assume that $\ell=0$, or $\ell>0$ and $\neq p$.
Let $G$ be a reductive group.
Let $\Hbu$ be the affine Hecke algebra of $G$ over $\k$.
Cherednik has introduced a double affine Hecke algebra $\Hb$, 
which may be viewed as an affine counterpart to $\Hbu$.
It is natural to guess that $\Hb$ takes some role 
in the representation theory of $\k G(\F)$.
More precisely, let $\Bc$ be the unipotent block in the category of smooth
representations of $\k G(\F)$,
i.e. the block containing the trivial representation.
We expect $\Bc$ to be equivalent
to some category of representations of $\Hb$. 

The main result of this paper is a step in this direction.
Assume that $G=\GL_n$.
Fix an Iwahori subgroup $I\subset G(\F)$.
Let $\Ic$ be the annihilator of the natural representation 
of the global Hecke algebra of $G(\F)$ in $\k[G(\F)/I]$.  
The full subcategory $\Bc'\subset\Bc$ consisting of representations
annihilated by $\Ic$ is an Abelian category.
Let $\Scbu$ be the quantized affine Schur algebra of $G$ over $\k$.
Recall that $\Hbu$, $\Scbu$ 
are algebras over the ring $\k[\zeta^{\pm 1}]$,
with $\zeta$ the quantum parameter,
while $\Hb$ is an algebra over $\k[\tau^{\pm 1},\zeta^{\pm 1}]$.
It is proved in \cite{Vi} that $\Bc'$ is equivalent to 
$(\Scbu|_{\zeta=q})$-$\Mc od$.
Note that $q$ is a root of unity in $\k^\times$ if $\ell>0$.
Given $h_0,u_0\in\k$, set $\tau_0=\e^{u_0}$, $\zeta_0=\e^{u_0h_0}$ 
and let $\Oc_{\zeta_0,\tau_0}\subset(\Hb|_{\tau=\tau_0,\zeta=\zeta_0})$-$mod$ 
be the `category $\Oc$'.
If $\k=\CC$ and $u_0$ is generic we prove an equivalence 
between blocks of $\Oc_{\zeta_0,\tau_0}$ 
and blocks of $(\Scbu|_{\zeta=\e^{h_0}})$-$mof$.
The whole categories are not equivalent because they have different centers.
To get the same categories one must
either replace $\Scbu$ by some elliptic analogue (which is not known so far),
or replace $\Oc_{\zeta_0,\tau_0}$ by the corresponding category 
$\Oc'_{h_0}$ of modules of the double affine graded Hecke algebra $\Hb'$.
See Section 5 for precise statements.
We conjecture that our equivalence 
is still true if $\k$ is any algebraically closed
field of characteristic $\ell>0$.

Roughly, the proof is as follows. 
We split $\Oc_{\zeta_0,\tau_0}$ as a direct sum of subcategories
$\Oc_{\zeta_0,\tau_0}=\Oplus_\ell{}^{\{\ell\}}\Oc_{\zeta_0,\tau_0}.$ 
Each summand is equivalent to a category of modules,
say ${}^{\{\l\}}\Oc_{h_0}$, of $\Hb'$.
The category ${}^{\{\l\}}\Oc_{h_0}$ 
is the limit of an inductive system of subcategories
${}^{\l}\Oc_{n,h_0}$ with $n\in\ZZ_{\geq 0}$.
Althought ${}^{\{\l\}}\Oc_{h_0}$ do not have enough projective objects,
the categories ${}^{\l}\Oc_{n,h_0}$ 
are generated by a family of projective modules
which are easily described.
We construct an exact functor 
$\Mc\,:\,{}^\l\Oc_{n,h_0}\to\Hbu_\CC$-$mof$
which is faithful on projective objects, under a mild restriction,
using the trigonometric Knizhnik-Zamolodchikov connection.
This functor is inspired from \cite{GGOR}.
In general we do not know how to compute the image by $\Mc$ 
of any projective generator.
However, in some particular cases including the type $A$ case,
this can be done via some deformation argument.

We may have proved our equivalence of categories
with the geometric technics used in \cite{V2}.
From this viewpoint, one is essentially reduced to prove the 
injectivity conjectured in \cite{V2, Remark 4.9}.
By {\it loc. cit.}, in type $A$,
the simple object in $\Oc_\CC$ 
are labelled by representations of a cyclic quiver,
and the Jordan-Holder multiplicities  of induced modules
are the value at one of certain Kazhdan-Lusztig polynomials of type $A^{(1)}$.
Our equivalence of categories may be viewed as an extension of these results.
However, the present approach is more powerful in the sense that  
the $K$-theoretic construction does not adapt easily to the case
of double affine Hecke algebras with several parameters.

\head 1. Notations\endhead
\subhead 1.1. Reminder on modules and categories\endsubhead
Let $\k$ be a principal domain of characteristic zero. 
We will mainly assume that $\k=\A,\F$ or $\CC$,
where $\A=\CC[[\varpi]]$ and $\F=\CC((\varpi))$.
Let $\k^\times\subset\k$ be the multiplicative group. 
Given a $\k$-algebra $\Ab$,
let $\Ab$-$\Modc$ be the category of left 
$\Ab$-modules which are free over $\k$,
$\Ab$-$mod$ be the full subcategory consisting of finitely generated modules,
$\Ab$-$mof$ be the full subcategory 
consisting of the modules of finite type over $\k$.

Given an Abelian category $\Ac$ and a 
full subcategory $\Nc\subset\Ac$ stable under subquotients
and extensions, let $\Ac/\Nc$ be the Serre quotient, see \cite{G}.
The category $\Ac/\Nc$ is Abelian and the obvious functor $Q\,:\,\Ac\to\Ac/\Nc$
is exact. Given an Abelian category $\Bc$ and an exact functor $F\,:\,\Ac\to\Bc$
such that $FM\simeq 0$ for all $M\in\Nc$, there is a unique exact functor
$G\,:\,\Ac/\Nc\to\Bc$ such that $F=G\circ Q$.
Conversely, given an Abelian category $\Cc$,
an exact functor $Q\,:\,\Ac\to\Cc$ is called a quotient functor
if and only if it induces an equivalence $\Ac/\ker Q\to\Cc$.
Clearly, $Q$ is a quotient functor if and only if for any
exact functor $F\,:\,\Ac\to\Bc$ such that $FM\simeq 0$ whenever $QM\simeq 0$
there is a unique exact functor $G\,:\,\Cc\to\Bc$ such that $F=G\circ Q$.
If $\Ac$ is Artinian and $P\in\Ac$ is projective, the functor
$\Hom_\Ac(P,\cdot)$ is a quotient functor from
$\Ac$ to the category of right $\End_\Ac P$-modules of finite length. 

\subhead 1.2. Reminder on root systems\endsubhead
Let $\Delta$ be an irreducible root system.
Let $\Delta_+\subset\Delta$ be a system of positive roots,
and $\Pi=\{\a_i\,;\,i\in I\}\subset\Delta_+$ be the simple roots. 
Let $\theta\in\Delta_+$ be the maximal root,
and $\rho={1\over 2}\sum_{\b\in\Delta_+}\b$.
The set of simple affine roots is 
$\{\a_i\,;\,i\in\hat I\}$, where $\hat I=I\cup\{\heartsuit\}$.
For any subset $J\subseteq\hat I$ set $\Pi_J=\{\a_i\,;\,i\in J\}$, 
$\Delta_J=\Delta\cap\ZZ\Pi_J$,
and $\Delta_{J,+}=\Delta_+\cap\Delta_J$.
Let $\Delta^\vee$, $\Delta_+^\vee$, etc., denote the corresponding
sets of coroots.
Let $\hat\Delta_\re=\Delta\times\ZZ$,
$\hat\Delta_\re^\vee=\Delta^\vee\times\ZZ$ 
be the sets of affine real roots and coroots.

Denote by $Y$, $Y^\vee$ the root and the coroot lattices,
by $X$, $X^\vee$ the weights and the coweights lattices.
Let $Y_+\subset Y$ be the semigroup generated by $\Delta_+$,
and write $Y_{\RR,+}$ for $\RR_{\geq 0}\otimes Y_+$.

Let $W$, $\hat W$ be the Weyl group and the affine Weyl group.
Let $s_\b\in W$ (resp. $s_{\hat\b}\in\hat W$)
be the reflection relatively to the root 
$\b\in\Delta$ (resp. $\hat\b\in\hat\Delta_\re$).
We write $s_i$ for $s_{\a_i}$. 
Recall that $\hat W=Y\rtimes W$. 
We write $x_\b$ for $(\b,0)$,
and $s_\heartsuit$ for $x_{\theta}s_\theta$. 
Let $\ell\,:\,W\to\ZZ_{\geq 0}$ be the length.
We write $\geq$ for the Bruhat order on $\hat W$.

Let $\Omega\subset\Aut(\hat W)$ be the group of diagram automorphisms,
and $\tilde W=\hat W\rtimes\Omega$ be the extended affine Weyl group.
Then $\tilde W\simeq X\rtimes W$. 
As above we write $x_\mu$ for $(\mu,0)$.
For each $\pi\in\Omega\setminus\{1\}$ 
let $\a_\pi\in\Pi$ be such that $\pi(s_\heartsuit)=s_{\a_\pi}$.
Let $\omega^\vee_\pi$ be the fundamental coweight dual to $\a_\pi$,
and $\omega_\pi$ be the corresponding fundamental weight.
We write $\pi$ for the element $(0,\pi)\in\tilde W$.
Let $w_\pi\in W$ be such that 
$\pi\in\tilde W$ is identified with $x_{\o_\pi}w_\pi\in X\rtimes W$.
Hence $w_\pi\theta=-\a_\pi$, and $w_\pi\a_i=\a_j$ if $i,j\neq\heartsuit$ 
and $\pi(s_i)=s_j$.

Set $X_\k=\k\otimes X$, $X^\vee_\k=\k\otimes X^\vee$. 
The group $W$ acts on $X_\k$ and $X^\vee_\k$ by
$s_\b\l=\l-(\l:\b^\vee)\b$ and
$s_\b\l^\vee=\l^\vee-(\b:\l^\vee)\b^\vee,$
where $(\,:\,)$ is the unique $\k$-linear
pairing $X_\k\times X^\vee_\k\to\k$ such that 
$(\a_i\,:\,\o_j^\vee)=\delta_{ij}$.
We write $\geq$ for the order on $X_\k$ such that
$\mu\geq\nu$ if and only if $\mu-\nu\in Y_+$.

Let $\Sb'$ be the symmetric algebra of $X^\vee_\k$.
Given $\l^\vee\in X^\vee$ set $\xi_{\l^\vee}=1\otimes \l^\vee\in\Sb'$.
Put $\xi_i=\xi_{\o_i^\vee}$,
and $\xi_{\a_\heartsuit^\vee}=1-\xi_{\theta^\vee}$.
The group $\tilde W$ acts on the $\k$-algebra $\Sb'$ by 
${}^{x_\mu w}\xi_{\l^\vee}=\xi_{w\l^\vee}-(\mu\,:\,w\l^\vee).$
The dual action on $X_\k$ is $x_\mu w(\l)=\mu+w\l.$
Put $\xi_{\hat\b^\vee}=\xi_{\b^\vee}+r$ if 
$\hat\b^\vee=(\b^\vee,r)\in\hat\Delta^\vee_\re$.
There is a unique $\tilde W$-action on $\hat\Delta_\re^\vee$
such that ${}^w\xi_{\hat\b^\vee}=\xi_{w\hat\b^\vee}$.

Set $T=\k^\times\otimes X^\vee$, and $T^\vee=\k^\times\otimes Y$.
In the following $\otimes$ means $\otimes_\ZZ$,
and $\e^z$ means $\exp(2i\pi z)$. 
For any $\l\in X_\k$ (resp. $\l^\vee\in X_\k^\vee$)
we write $\l_j$ for the element $(\l\,:\,\o^\vee_j)$ in $\k$
(resp. $\l^\vee_j$ for $(\o_j\,:\,\l^\vee)$)
and $\e^\l$ for the element $\prod_j\e^{\l_j}\otimes\a_j$ in $T^\vee$
(resp. $\e^{\l^\vee}$ for $\prod_j\e^{\l^\vee_j}\otimes\a^\vee_j\in T$).

Fix $\tau\in\k^\times$.
Fix a $m$-th root $\tau^{1/m}$ of $\tau$, 
where $m$ is the least natural number such that $(X\,:\,X^\vee)=(1/m)\ZZ$.

Set $\Sb=\k X^\vee$.
Given $\l^\vee\in X^\vee$ let $y_{\l^\vee}$ be the
corresponding element in $\Sb$.
Put $y_i=y_{\o_i^\vee}$,
and $y_{\a_\heartsuit^\vee}=\tau y_{-\theta^\vee}$.
There is a unique ring isomorphism $\Sb\simeq\k[T^\vee]$ 
taking $y_{\l^\vee}$ to the function 
$z\otimes\g\mapsto z^{(\g\,:\,\l^\vee)}$.
The group $\tilde W$ acts on the $\k$-algebra $\Sb$ by 
${}^{x_\mu w}y_{\l^\vee}=y_{w\l^\vee}\tau^{-(\mu\,:\,w\l^\vee)}.$
The dual action on $T^\vee$ is 
$x_\mu w(z\otimes\g)=(z\otimes w\g)(\tau\otimes\mu).$
Put $y_{\hat\b^\vee}=y_{\b^\vee}\tau^r$ if 
$\hat\b^\vee=(\b^\vee,r)\in\hat\Delta^\vee_\re$.

Set $\Rb=\k Y$.
Given $\b\in Y$ let $x_\b$ denote also the
corresponding element in $\Rb$.
We write $x_i$ for $x_{\a_i}$.
There is a unique ring isomorphism $\Rb\simeq\k[T]$ 
taking $x_\b$ to the function 
$z\otimes\l^\vee\mapsto z^{(\b\,:\,\l^\vee)}$.
The group $W$ acts on the $\k$-algebra $\Rb$ by ${}^w x_\b=x_{w\b}$.
Let $X^\vee_\k\times\Rb\to\Rb$, $(\xi,f)\mapsto\partial_\xi f$
be the unique $\k$-linear action such that
$\partial_{\xi_{\l^\vee}}(x_\b)=(\b\,:\,\l^\vee)x_\b$. 

Given a root $\b$ let $\vartheta_\b\,:\,\Rb\to\Rb$,
$\vartheta_{\b^\vee}\,:\,\Sb'\to\Sb'$ be the $\k$-linear operators such that
$$\vartheta_{\b^\vee}(p)={p-{}^{s_\b}p\over\xi_{\b^\vee}},\quad
\vartheta_\b(f)={f-{}^{s_\b}f\over 1-x_{-\b}}.$$

\subhead 1.3. Reminder on affine Weyl groups\endsubhead
For each subset $J\subsetneq\hat I$ 
let $W_J\subset\hat W$ be the subgroup generated by $\{s_i\,;\,i\in J\}$.
It is finite.
Let $W^J\subseteq W$ be the set of elements $v$ such that
$\ell(vu)=\ell(v)+\ell(u)$ for each $u\in W_J$.

If $\ell\in T^\vee$ we put $\hat W_\ell=\{w\in\hat W\,;\,w\ell=\ell\}$
and $W_\ell=W\cap\hat W_\ell$.
If $\l\in X_\k$ we put
$\hat W_\l=\{w\in\hat W\,;\,w\l=\l\}$ and $W_\l=W\cap\hat W_\l$.
The group $\hat W_\l$ is finite.
If $\tau$ is not a root of unity then $\hat W_\ell$ is also finite.
Let $\hat n_\l,n_\ell$ be the number of elements in 
$\hat W_\l$, $W_\ell$ respectively.

If $\k$ is an algebraically closed field the group $W_\ell$
is generated by reflections, see \cite{SS, Chap. II, Theorem $\ $4.2}.

\proclaim{Lemma}
(i)
Any finite subgroup in $\hat W$ is conjugate into some $W_J$. 

(ii)
$W_\l=W_{\e^\l}\iff W_\l=\hat W_\l$.
\endproclaim

\noindent{\sl Proof.}
Given a finite subgroup $W'\subset\hat W$, 
and $x$ an element in the interior of the Tits cone,
the stabilizer in $W$ of $\sum_{w\in W'}wx$ contains $W'$,
and is $\hat W$-conjugate into $W_J$ for some $J$ by 
\cite{K, Proposition 3.12}.  
Claim $(ii)$ is obvious because $W_{\e^\l}=\{w\in W\,;\,\l-w\l\in Y\}$,
and $\hat W_\l=\{x_{\l-w\l}w\,;\,w\in W_{\e^\l}\}$.
\qed

\head 2. The category $\Oc'$\endhead
\subhead 2.1. The category $\Oc'$\endsubhead
Fix $h_{\hat\b}\in\k$ for each $\hat\b\in\hat\Delta_\re$, 
such that $h_{\hat\b}=h_{\a_i}$ if 
$\hat\b\in\tilde W\a_i$ and $i\in\hat I$. 
We write $h_i$ for $h_{\a_i}$.
Let $\Hb'$ be the degenerate double affine Hecke algebra.
Recall that $\Hb'$ is the $\k$-algebra generated by $\k\hat W$ and $\Sb'$
with the relations
$$s_i p-{}^{s_i}p\,s_i=h_i\vartheta_{\a_i^\vee}(p),
\leqno(2.1.1)$$
for all $i\in\hat I$, $p\in\Sb'$. 
Then
$$\xi f-f\xi=\partial_\xi(f)-\sum_{\b\in\Delta_+}h_\b(\b\,:\,\xi)
\vartheta_\b(f)s_\b,\qquad
\forall f\in\Rb,\ \forall\xi\in X_\k^\vee.
\leqno(2.1.2)$$
The product yields an isomorphism
$\Rb\otimes_\k\k W\otimes_\k\Sb'\to\Hb'$.
There is a unique action of $\Omega$ on $\Hb'$ by algebra automorphisms
such that $\pi\in\Omega$ acts on $\hat W$ and $\Sb'$ as in 1.2.

Let $\Oc'\subset\Hb'$-$mod$ be the full subcategory consisting of
the modules which are locally finite respectively to $\Sb'$.
To avoid some ambiguity we may write $\Hb'_\k$ for $\Hb'$.

Set $\la\mu\ra=\{p-p(\mu)\,;\,p\in\Sb'\}$, 
and $\la E\ra=\bigcap_{\mu\in E}\la\mu\ra$
for each finite subset $E\subset\hat W\l$.
Let ${}^{\{\l\}}\Oc\subset\Oc'$ be the full subcategory consisting of
the modules $M$ such that for each element $m\in M$ there is a finite subset
$E\subset\hat W\l$ and an integer $n>0$ with 
$\la E\ra^n m=\{0\}$.

\proclaim{Proposition}
(i)
$\Oc'\subset\Hb'$-$\Modc$ is a Serre subcategory.

(ii)
If $\k$ is an algebraically closed field
then $\Oc'=\bigoplus_\l{}^{\{\l\}}\Oc,$
where $\l$ varies in a set of representatives of the $\hat W$-orbits in $X_\k$.
\endproclaim

\noindent{\sl Proof.}
Any object in $\Oc'$ is finitely generated over $\Rb$,
because $W$ is finite, $\Hb'=\Rb\cdot\k W\cdot\Sb'$ and 
a module in $\Oc'$ is finitely generated over $\Hb'$ 
and locally finite over $\Sb'$.
Hence the category $\Oc'$ is Abelian, because $\Rb$ is a Noetherian ring.
The category $\Oc'$ is obviously closed by subquotients and extensions.

For each module $M$ in $\Oc'$ let ${}^{\{\l\}}M\subset M$ be the subspace  
consisting of the elements $m\in M$ such that there is a finite subset
$E\subset\hat W\l$ and an integer $n>0$ with $\la E\ra^n m=\{0\}$.
Clearly ${}^{\{\l\}}M$ lies in ${}^{\{\l\}}\Oc$.
We have $M=\sum_\l{}^{\{\l\}}M$ because the 
$\Sb'$-action on $M$ is locally finite, 
and this sum is obviously direct.
Claim $(ii)$ follows. 
\qed

\vskip3mm

For any group $G$ acting linearly on $X_\k$ and any $\l\in X_\k$,
let $[\l]_{G,\k}\subset\Sb'$ be the ideal generated by
$\la\l\ra^G$ (=the $G$-invariant elements in $\la\l\ra$).
We write $[\l]$ (or $[\l]_\k$ if necessary)
for $[\l]_{\hat W_\l,\k}$.
Set $\Sb_\l=\Sb'/[\l]$ (or $\Sb_{\l,\k}$ if necessary).
Note that $\Sb_\l$ is of finite type over $\k$ because $\hat W_\l$ is finite,
see \cite{B, chap. V, \S 1, $n^o$ 9, Th\'eor\`eme 2}.
If $\k$ is a local ring then $\Sb_\l$ is also a local ring.
In this case let $\mb_\l\subset\Sb_\l$ be the maximal ideal.

Let $E\subset\hat W\l$ be finite.
Set $[E]=\bigcap_{\mu\in E}[\mu]$.
The quotient $\Sb_E=\Sb'/[E]$ is of finite type over $\k$.
If $\Sb_\l$ is free over $\k$ then $\Sb_E$ is also free because
it embeds in $\bigoplus_{\mu\in E}\Sb_\mu$ and $\k$ is principal.
If $\k$ is a field then $\Sb_E=\bigoplus_{\mu\in E}\Sb_\mu$.
If confusion is unlikely from the context 
we write $p$ again for the image in $\Sb_E$ of an element $p\in\Sb'$. 

Let ${}^\l\Oc\subset{}^{\{\l\}}\Oc$ be the full subcategory consisting of
the modules such that 
for each element $m$ there is a finite subset 
$E\subset\hat W\l$ with $[E] m=\{0\}$.
For a future use we prove the following technical lemma.

\proclaim{Lemma}
If $\Sb_\l$ is torsion free over $\k$ there is a finite subset
$F\subset\hat W\l$ containing $E$ such that
$[F] s_i\subseteq s_i[E]+[E]$
in $\Hb'$.
\endproclaim

\noindent{\sl Proof.}
Fix a finite subset $F\subset\hat W\l$ containing $E$ such that $s_iF=F$.
We prove that $[F] s_i\subseteq s_i[F]+[F].$
For each $p\in[F]$ we have $p\,s_i=s_i\,{}^{s_i}\!p+\vartheta_{\a_i^\vee}(p)$
by (2.1.1).
Hence we must prove that $\vartheta_{\a_i^\vee}([F])\subseteq [F]$, i.e. that
$\vartheta_{\a_i^\vee}([\mu]\cap[s_i\mu])\subseteq [\mu]\cap [s_i\mu]$ for each
$\mu\in F$.

If $s_i\mu=\mu$ 
we are done because $[\mu]$ is generated by $\la\mu\ra^{\hat W_\mu}$,
for all $p_1,p_2\in\Sb'$ we have
$\vartheta_{\a_i^\vee}(p_1p_2)=\vartheta_{\a_i^\vee}(p_1)p_2+
{}^{s_i}p_1\vartheta_{\a_i^\vee}(p_2)$, and
$\vartheta_{\a_i^\vee}(\la\mu\ra^{\hat W_\mu})=\{0\}$.

Assume that $s_i\mu\neq\mu$. 
Fix $p\in[\mu]\cap[s_i\mu]$. 
It suffices to prove that $\vartheta_{\a_i^\vee}(p)\in [\mu]$.
We have $\xi_{\a_i^\vee}\vartheta_{\a_i^\vee}(p)=0$ in $\Sb_{\mu,\k}$.
Let $\K$ be the fraction field of $\k$.
Then $\xi_{\a_i^\vee}$ is invertible in $\Sb_{\mu,\K}$ because
$\Sb_{\mu,\K}$ is a local ring and $\xi_{\a_i^\vee}\notin\la\mu\ra$.
Hence $\vartheta_{\a_i^\vee}(p)=0$ in $\Sb_{\mu,\k}$,
because $\Sb_{\mu,\k}$ is torsion free over $\k$.
\qed

\vskip3mm

\noindent{\bf Remark.}
Let $G$ be a finite group acting linearly on $X_\CC$.
Fix $\l\in X_\A$ whose image, $\l_0$, in $X_\CC$ is fixed by $G$. 
Then the algebra $\Rb=\Sb'_\A/[\l]_{G,\A}$ is free over $\A$,
and $\Rb\otimes_\A\CC=\Sb'_\CC/[\l_0]_{G,\CC}$,
because the graded ring associated to the decreasing filtration 
$(\Rb\varpi^n)$ of $\Rb$ is isomorphic to 
$(\Sb'_\CC/[\l_0]_{G,\CC})\otimes_\CC\A$.
Indeed, for each $n$, the obvious map
$\Sb'_\A\varpi^n\to\Sb'_\CC\varpi^n$
takes $\Sb'_\A\varpi^{n+1}+[\l]_{G,\A}\varpi^n$
into $[\l_0]_{G,\CC}\varpi^n$,
and the resulting ring homomorphism
$$\Sb'_\A\varpi^n/(\Sb'_\A\varpi^{n+1}+[\l]_{G,\A}\varpi^n)
\to(\Sb'_\CC\varpi^n)/([\l_0]_{G,\CC}\varpi^n)$$
is invertible (use averages over $G$).

\subhead 2.2. Projective modules in $\Oc'$\endsubhead
For each $\mu\in\hat W\l$ we set $P(\mu)=\Hb'/\Hb'[\mu]$.
To avoid some confusion we may write $P(\mu)_\k$ for $P(\mu)$.
Let $1_\mu\in P(\mu)$ be the image of the unity 
by the obvious projection $\Hb'\to P(\mu)$.

For a future use we set $M_\mu=\{m\in M\,;\,[\mu]m=0\}$
for each $\Hb'$-module $M$.
If $\k$ is a field and $M$ lies in ${}^\l\Oc$ then 
$M=\bigoplus_{\mu\in\hat W\l}M_\mu$,
because for any $m\in M$ the map $\Sb'\to M$, $p\mapsto pm$ factors through
$\Sb_E\to M$ for a finite set $E\subset\hat W\l$,
and $\Sb_E=\Oplus_{\mu\in E}\Sb_\mu$.

\proclaim{Proposition}
Assume that $\k$ is a field.

(i)
$P(\mu)$ is a projective object in ${}^\l\Oc$.

(ii)
The category ${}^\l\Oc$ is generated by the modules
$P(\mu)$ with $\mu\in\hat W\l$.

(iii) 
The category $\Oc'$ is Artinian, 
and there are a finite number of simple objects in ${}^\l\Oc$.
\endproclaim

\noindent{\sl Proof.}
For each $w\in\hat W$ there is a finite subset $E\subset\hat W\l$ such that
$[E] w\subset\sum_{w'\leq w}w'[\mu]$ by Lemma 2.1. 
Then, $[E]w\,1_\mu=0$ because $[\mu]1_\mu=0$.
Therefore $P(\mu)$ belongs to ${}^\l\Oc$.

Given a map $f\,:\,M\to N$ in ${}^\l\Oc$,
we have $f(M)=\Oplus_\mu f(M_\mu)=\Oplus_\mu f(M)_\mu$,
and $f(M_\mu)\subseteq f(M)_\mu$ for each $\mu$.
Therefore $f(M_\mu)=f(M)_\mu$.
Thus $P(\mu)$ is projective because $M_\mu=\Hom_{\Hb'}(P(\mu),M)$ for each $M$.
Claim $(i)$ is proved.

Each object $M$ in ${}^\l\Oc$ is a quotient 
of a direct sum of modules isomorphic to some $P(\mu)$,
because $M=\bigoplus_{\mu\in\hat W\l}M_\mu$.
Since $M$ is finitely generated it is indeed the quotient 
of a finite direct sum of these modules.
Claim $(ii)$ is proved.

To prove that $\Oc'$ is Artinian 
it is sufficient to check that $P(\mu)$ has finite length over $\Hb'$
for each $\mu$. We have
$${}^w\!\xi\, w-w\,\xi\in\Sum_{w'<w}\k w',\quad
\forall\xi\in X_\k^\vee,\,w\in\hat W.\leqno(2.2.1)$$
Let $P(\mu)_{\leq w}\subseteq P(\mu)$ be the right $\Sb_\mu$-submodule 
spanned by $\{w'1_\mu\,;\,w'\leq w\}$.
Then $(P(\mu)_{\leq w})$ is a filtration of $P(\mu)$
by left $\Sb'$-submodule, by (2.2.1).
Let $P(\mu)_\bullet$ be the associated graded.
The $\Hb'$-action on $P(\mu)$ yields a $\k\hat W\rtimes\Sb'$-action on
$P(\l)_\bullet$ by (2.2.1), where $\k\hat W\rtimes\Sb'$ is the semi-direct 
product relative to the $\k\hat W$-action on $\Sb'$ in 1.2.
It is sufficient to prove that $P(\l)_\bullet$ has a finite length
over $\k\hat W\rtimes\Sb'$.
Note that $P(\l)_\bullet$ is isomorphic to
$\bigoplus_{\mu\in\hat W\l}(\Sb_\mu)^{\oplus\hat n_\mu}$ over $\Sb'$,
and that a $(\k\hat W\rtimes\Sb')$-submodule of $P(\l)_\bullet$ 
is a sum of $\Sb'$-submodules 
$U_\mu\subseteq(\Sb_\mu)^{\oplus\hat n_\mu}$ 
such that $w(U_\mu)=U_{w\mu}$ for all $w\in\hat W$.
Thus the length of $P(\l)_\bullet$ is bounded 
by the length of $(\Sb_\l)^{\oplus\hat n_\l}$ over $\Sb'$.
Hence it is finite because $\k$ is a field.

By $(ii)$, the last part of $(iii)$ is a consequence of Proposition 2.3 below.
\qed

\vskip3mm

\noindent{\bf Remarks.}
$(i)$
If $\k$ is a field, 
simple objects in ${}^\l\Oc$ have projective covers.
However they do not have finite projective resolutions in general.

$(ii)$
In general $P(\mu)$ is not indecomposable over $\Hb'$.

$(iii)$
Assume that $\k$ is a field.
Each object $M\in{}^{\{\l\}}\Oc$ 
has a filtration whose associated graded lies in ${}^{\l}\Oc$
(consider the submodule $\{m\in M\,;\,\exists E\,\roman{s.t.}\,[E]m=0\}$,
which lies in ${}^{\l}\Oc$,
and use the fact that $M$ has a finite length).
If $\k$ is algebraically closed and $M\in\Oc'$ is simple 
then it lies in ${}^\l\Oc$ for some $\l\in X_\k$,
because it lies in ${}^{\{\l\}}\Oc$ for some $\l\in X_\k$,
hence it has a filtration whose associated graded lies in ${}^{\l}\Oc$.

\subhead 2.3. Intertwiners in $\Oc'$\endsubhead
Assume that $\k$ is a field.
For any reduced decomposition 
$w=s_{i_1}s_{i_2}\cdots s_{i_r}\in\hat W$ set
$\phi'_w=\phi'_{i_1}\cdots\phi'_{i_{r-1}}\phi'_{i_r}\in\Hb',$
with $\phi'_i=s_i\xi_{\a_i^\vee}-h_i$ for all $i\in\hat I$.
Recall that ${}^w p\,\phi'_w=\phi'_wp$ for all $p\in\Sb'.$
The intertwining operator 
$\Phi'_w(\mu)\,:\,P(w\mu)\to P(\mu)$ is the unique 
$\Hb'$-homomorphism taking
$1_{w\mu}$ to $\phi'_w 1_\mu$.

\proclaim{Lemma}
The operator $\Phi'_{s_i}(\mu)$ is invertible if and only if
$(\mu\,:\,\a_i^\vee)\neq\pm h_i$.
\endproclaim

\noindent{\sl Proof.}
Set $\psi'_i(\mu)=(\mu\,:\,\a_i^\vee)s_i-h_i$.
Let
$\Psi'_{s_i}(\mu)\,:\,P(s_i\mu)\to P(\mu)$ be the unique 
$\k\hat W$-homomorphism taking $wp$ to 
$w\psi'_i(\mu)p$ for each $w\in\hat W$.
The $\k$-modules $P(\mu)^{\geq k}=\k\hat W\otimes_\k(\mb_\mu)^k$, 
with $k\geq 0$,
form a finite decreasing filtration of $P(\mu)$.
We have $\Psi'_{s_i}(\mu)(P(s_i\mu)^{\geq k})\subseteq P(\mu)^{\geq k}$ 
for each $k$,
and $\Psi'_{s_i}(\mu)$, $\Phi'_{s_i}(\mu)$ 
coincide in the associated graded spaces.
Hence $\Phi'_{s_i}(\mu)$ is invertible 
if and only if $\Psi'_{s_i}(\mu)$ is invertible.
The lemma follows.
\qed

\vskip3mm

For each  $\hat\b^\vee\in\hat\Delta^\vee_\re$ we put
$H_{\hat\b^\vee}=\{\mu\in X_\RR\,;\,\xi_{\hat\b^\vee}(\mu)=0\}$.
The connected components of
$X_\RR\setminus\bigcup_{\hat\b^\vee\in\hat\Delta_\re^\vee}H_{\hat\b^\vee}$
are the alcoves.
Let $A_+$ be the alcove containing $\rho/k$ if $k\gg 1$,
and $A_w=\{w^{-1}\mu\,;\,\mu\in A_+\}$ for each $w\in\hat W$.

The set $\Hc_\l=\{\hat\b^\vee\in\hat\Delta^\vee_\re\,;
\,\xi_{\hat\b^\vee}(\l)=\pm h_{\hat\b}\}$ is finite.
Set 
$U_\l=X_\RR\setminus\bigcup_{\hat\b^\vee\in\Hc_\l}H_{\hat\b^\vee}$.
The group $\hat W_\l$ acts on $U_\l$. 
An affine domain is a minimal subset in $U_\l$ containing
a connected component and stable by $\hat W_\l$. 
Let $D_w$ be the unique affine domain containing $A_w$,
and let $\Dc$ be the set of affine domains.

\proclaim{Proposition}
(i)
The $\Hb'$-modules $P(w_1\l)$, $P(w_2\l)$ are isomorphic 
if $D_{w_1}=D_{w_2}$.

(ii)
The modules 
$P(\l)$, $P(w\l)$ have the same composition factors
for all $w\in\hat W$.
\endproclaim

\noindent{\sl Proof.}
Fix $w\in\hat W$ and $i\in\hat I$.
The intertwining operator 
$\Phi'_{s_i}(w\l)\,:\,P(s_iw\l)\to P(w\l)$ 
is invertible if and only if 
$\xi_{w^{-1}\a_i^\vee}(\l)\neq\pm h_i$.
Thus $\Phi'_{w_1w_2^{-1}}(w_2\l)\,:\,P(w_1\l)\to P(w_2\l)$ 
is invertible if $A_{w_1v_1}$, $A_{w_2v_1}$ are in the same affine domain
for some $v_1,v_2\in\hat W_\l$.
This gives $(i)$.

Fix $\l$, $w$. 
The modules $P(\l)$, $P(w\l)$ are isomorphic 
for generic parameters $h_i$ by $(i)$.
Hence $(ii)$ follows by a standard argument,
see \cite {CG, Lemma 2.3.4} for instance. 
\qed

\subhead 2.4. Induction\endsubhead
For any subset $J\subsetneq\hat I$,
the $\k$-submodule $\Hb'_J=\k W_J\cdot\Sb'\subset\Hb'$ is a subring by (2.2.1).
Set $\Hbu'=\Hb'_I$ and $\Ocu'=\Hbu'$-$\mof$.
For each $\Hbu'$-module $M$ let $\Ic(M)=\Hb'\otimes_{\Hbu'}M.$ 
Put $\Pu(\l)=\Hbu'/\Hbu'[\l]$. Then $\Ic(\Pu(\l))=P(\l)$.
If $M$ is finitely generated over $\Hbu'$ then 
$\Ic(M)$ is finitely generated over $\Hb'$.
If $M$ is locally finite over $\Sb'$ then 
$\Ic(M)$ is also locally finite over $\Sb'$ by (2.2.1), because
$\Ic(M)\simeq\k\hat W\otimes_{\k W}M$.
Thus $\Ic$ factors through a functor $\Ocu'\to\Oc'$.

\subhead 2.5. The category $\Oc$\endsubhead
We do not assume anymore that $\k$ is a field.
Fix $\zeta_{\hat\b}\in\k^\times$ for each $\hat\b\in\hat\Delta_\re$, 
such that $\zeta_{\hat\b}=\zeta_{\a_i}$ if 
$\hat\b\in\tilde W\a_i$ and $i\in\hat I$. 
We write $\zeta_i$ for $\zeta_{\a_i}$.
Let $\Hb$ be the corresponding double affine Hecke algebra.
It is the k-algebra generated by $\Sb$
and the elements $t_w$ with $w\in\hat W$,
modulo the following relations
$$(t_i-\zeta_i)(t_i+1)=0,\quad t_vt_w=t_{vw},$$
$$t_iy_{\l^\vee}-{}^{s_i}y_{\l^\vee}t_i=
(\zeta_i-1)(y_{\l^\vee}-{}^{s_i}y_{\l^\vee})(1-y_{-\a_i^\vee})^{-1},$$
if $\ell(vw)=\ell(v)+\ell(w)$ and $t_i=t_{s_i}$.
We may write $\Hb_\k$ for $\Hb$ if necessary.
There is a unique action of $\Omega$ on $\Hb$ by algebra automorphisms
such that $\pi(t_w)=t_{\pi(w)}$ for each $w\in\hat W$, and
$\pi(p)={}^\pi p$ for each $p\in\Sb$.

For any reduced decomposition 
$w=s_{i_1}s_{i_2}\cdots s_{i_r}\in\hat W$ set
$\phi_w=\phi_{i_1}\phi_{i_2}\cdots\phi_{i_r}\in\Hb,$
with  $\phi_i=t_i(y_{-\a_i^\vee}-1)+\zeta_i-1$ for all $i\in\hat I$.
Recall that
${}^wp\phi_w=\phi_wp$
for all $p\in\Sb.$

Fix $\ell\in T^\vee$.
Let $\la\ell\ra=\{p-p(\ell)\,;\,p\in\Sb\}$.
For any group $G$ acting on $T^\vee$,
let $[\ell]_{G,\k}\subset\Sb$ be the ideal generated by $\la\ell\ra^G$.
We write $[\ell]$ (or $[\ell]_\k$ if necessary)
for $[\ell]_{\hat W_\ell,\k}$.
Set $\la E\ra=\bigcap_{m\in E}\la m\ra$, $[E]=\bigcap_{m\in E}[m]$ 
if $E\subset\hat W\ell$ is finite.

Let $\Oc\subset\Hb$-$mod$ be the full subcategory consisting of
the modules which are locally finite respectively to $\Sb$.
Let ${}^{\{\ell\}}\Oc\subset\Oc$ (resp. ${}^\ell\Oc\subset\Oc$)
be the full subcategory consisting of
the modules $M$ such that for each element $m\in M$ there is a finite subset
$E\subset\hat W\ell$ such that $\la E\ra^n m=\{0\}$ if $n\gg 0$ 
(resp. such that $[E]\,m=\{0\}$).

If $\k=\CC$ we write $h_{0i},\zeta_{0i},\tau_0$ for $h_i,\zeta_i,\tau$.
Given $u_0\in\CC^\times$ such that 
$$u_0(k+(\l_0\,:\,\b^\vee))\notin\ZZ\setminus\{0\},\quad
\forall\b^\vee\in\Delta^\vee\cup\{0\},\ \forall k\in\ZZ+\Sum_i\ZZ h_{0i},
\leqno(2.5.1)$$
we set $\zeta_{0i}=\e^{u_0h_{0i}},$ $\tau_0=\e^{u_0}$.
Let $\Gamma\subset\CC^\times$ 
be the subgroup generated by $\e^{u_0}$, $\e^{u_0h_{0i}}$.
Fix $\ell_0\in T^\vee.$
The set 
$\Delta_{(\ell_0)}^\vee=
\{\a^\vee\in\Delta^\vee\,;\,y_{\a^\vee}(\ell_0)\in\Gamma\}$ 
is a root system. 
Let $\Delta_{(\ell_0)}\subseteq\Delta$ be the dual root system.
Let $\hat W_{(\ell_0)}$ be the affine Weyl 
group associated to $\Delta_{(\ell_0)}$.
Let $\Hb'_{(\ell_0)}$ be the degenerated double affine
Hecke algebra generated by $\hat W_{(\ell_0)}$ and $\Sb'$,
modulo the relation analoguous to (2.1.1), 
relatively to the set of parameters 
$\{h_{\hat\b}\,;\,\hat\b\in\Delta_{(\ell_0)}\times\ZZ\}$. 
For any $\l_0\in X_\CC$ let the category 
${}^{\{\l_0\}}\Oc_{(\ell_0)}\subset\Hb'_{(\ell_0)}$-$mod$ 
be as in 2.1.
Fix $\l_0$ such that 
$y_{\a^\vee}(\ell_0)=\e^{u_0(\l_0:\a^\vee)}$
for each $\a^\vee\in\Delta_{(\ell_0)}^\vee$.

\proclaim{Proposition}
(i)
$\Oc$ is a Serre subcategory of $\Hb$-$\Modc$.

(ii)
If $\k$ is an algebraically closed field then
$\Oc=\bigoplus_\ell{}^{\{\ell\}}\Oc$,
where $\ell$ varies in a set of representatives 
of the $\hat W$-orbits in $T^\vee$.

(iii) 
Set $\k=\CC$.
Assume that $h_{0i},\zeta_{0i},\tau_0,\l_0$ are as above.
If $\hat W_{\ell_0}$ is generated by reflections
there is an equivalence of categories
${}^{\{\ell_0\}}\Oc\simeq {}^{\{\l_0\}}\Oc_{(\ell_0)}$.
Moreover, if $\Delta^\vee_{(\ell_0)}$ is also the set of coroots $\a^\vee$ 
such that $(\l_0\,:\,\a^\vee)\in\ZZ+\sum_i\ZZ h_{0i}$ then the categories
${}^{\{\l_0\}}\Oc_{(\ell_0)}$ and ${}^{\{\l_0\}}\Oc$ are equivalent.
\endproclaim

\noindent{\sl Proof.}
Claims $(i)$, $(ii)$ are proved as in 2.1-3. 
Claim $(iii)$ is `well-known', but there is no proof in the litterature.
It is proved as in \cite{L}, to which we refer for details. 
The proof consists of two parts 
(corresponding to the two reductions in \cite{L}),
the first of which being an isomorphism
between some completion
of $\Hb'$, $\Hb$ similar to Cherednik's isomorphism.

\vskip2mm

\noindent $(A)$
The rings $\Sb'/\langle E\rangle^n$, 
with $E\subset\hat W_{(\ell_0)}\cdot\lambda_0$ finite and $n\geq 0$,
form an inverse system.
Let ${}^{\{\l_0\}}\Sb_{(\ell_0)}$ be the projective limit.
Set also ${}^{\{\ell_0\}}\Sb_{(\ell_0)}=\pro\,\Sb/\langle E\rangle^n$,
with $E\subset\hat W_{(\ell_0)}\cdot\ell_0$ finite and $n\geq 0$.
We have 
$\hat W_{(\ell_0)}\cap\hat W_{\ell_0}=\hat W_{(\ell_0)}\cap\hat W_{\l_0}$,
because for any element $w\in\hat W_{(\ell_0)}$ we have 
$$\matrix
w\l_0=\l_0&\iff 
(w\l_0\,:\,\a^\vee)=(\l_0\,:\,\a^\vee),\ 
\forall\a^\vee\in\Delta^\vee_{(\ell_0)}
\hfill\cr
&\iff 
y_{\a^\vee}(w\ell_0)=y_{\a^\vee}(\ell_0),\ 
\forall\a^\vee\in\Delta^\vee_{(\ell_0)}
\hfill\cr
&\iff 
w\ell_0=\ell_0 
\hfill\cr
\endmatrix$$
since $u_0\notin\QQ$.
Hence there is a bijection 
$\hat W_{(\ell_0)}\cdot\ell_0\simeq\hat W_{(\ell_0)}\cdot\l_0$
which is compatible with the $\hat W_{(\ell_0)}$-actions. 
It yields a ring isomorphism 
${}^{\{\l_0\}}\Sb_{(\ell_0)}\simeq{}^{\{\ell_0\}}\Sb_{(\ell_0)}$ 
which is compatible
with the $\hat W_{(\ell_0)}$-actions.
Let $\Kb'$, ${}^{\{\l_0\}}\Kb_{(\ell_0)}$, $\Kb$, 
${}^{\{\ell_0\}}\Kb_{(\ell_0)}$ 
be the fraction fields of $\Sb'$, ${}^{\{\l_0\}}\Sb_{(\ell_0)}$, 
$\Sb$, ${}^{\{\ell_0\}}\Sb_{(\ell_0)}$.
Let $\Hb_{(\ell_0)}$ be the double affine Hecke algebra 
corresponding to $\Hb'_{(\ell_0)}$. 
Set 
$${}^{\{\l_0\}}\Hb_{(\ell_0)}=
{}^{\{\l_0\}}\Sb_{(\ell_0)}\otimes_{\Sb'}\Hb'_{(\ell_0)},\quad
{}^{\{\ell_0\}}\Hb_{(\ell_0)}=
{}^{\{\ell_0\}}\Sb_{(\ell_0)}\otimes_\Sb\Hb_{(\ell_0)}.$$
Set also
$${}^{\{\l_0\}}\Hb^\Kb_{(\ell_0)}=
{}^{\{\l_0\}}\Kb_{(\ell_0)}\otimes_{\Sb'}\Hb'_{(\ell_0)},\quad
{}^{\{\ell_0\}}\Hb^\Kb_{(\ell_0)}=
{}^{\{\ell_0\}}\Kb_{(\ell_0)}\otimes_\Sb\Hb_{(\ell_0)}.$$
For each $w\in\hat W$ let 
$\varphi'_w\in\Kb'\phi'_w$ 
(resp. $\varphi_w\in\Kb\phi_w$)
be normalized so that the map 
$w\mapsto\varphi'_w$ 
(resp. $w\mapsto\varphi_w$)
is a group homomorphism.
The intertwiner $\varphi_w$ is denoted by $G_w$ in \cite{C2}.
An element in 
${}^{\{\l_0\}}\Hb^\Kb_{(\ell_0)}$
is a finite sum $\sum_wp_w\varphi'_w$
with $w\in\hat W_{(\ell_0)}$ and $p_w\in{}^{\{\l_0\}}\Kb_{(\ell_0)}.$
By Lemma 2.1 there is a unique $\CC$-algebra structure on
${}^{\{\l_0\}}\Hb_{(\ell_0)}$, 
${}^{\{\l_0\}}\Hb^\Kb_{(\ell_0)}$ 
extending $\Hb'_{(\ell_0)}$.
Idem for ${}^{\{\ell_0\}}\Hb_{(\ell_0)}$, 
${}^{\{\ell_0\}}\Hb^\Kb_{(\ell_0)}$. 
There is a unique ring isomorphism
$${}^{\{\l_0\}}\Hb^\Kb_{(\ell_0)}\to
{}^{\{\ell_0\}}\Hb^\Kb_{(\ell_0)}
\ \roman{such\ that}\ 
\varphi'_w\mapsto\varphi_w,\ \forall w\in\hat W_{(\ell_0)}.$$
This isomorphism takes ${}^{\{\l_0\}}\Hb_{(\ell_0)}$ 
onto ${}^{\{\ell_0\}}\Hb_{(\ell_0)}$
(use (2.5.1) and see \cite{L, Theorem 9.3} for details).

\vskip2mm

\noindent $(B)$
Given $\ell\in\hat W\ell_0$, let $\Delta_{(\ell)}$ be the root system dual to
$\{\a^\vee\in\Delta^\vee\,;\,y_{\a^\vee}(\ell)\in\Gamma\}$.
Note that $\Delta_{(\ell)}=\Delta_{(\ell')}$ if $\ell,\ell'$ belong to the same 
$\Gamma\otimes Y$-coset.
The elements $\ell,\ell'$ are said to be equivalent
if they belong to the same $\Gamma\otimes Y$-coset 
and to the same $\hat W_{(\ell)}$-orbit.
Let $\Pc$ be the set of equivalence classes in $\hat W\ell_0$.
We write $(\ell)\in\Pc$ for the class of $\ell$, i.e.
$(\ell)=\hat W_{(\ell)}\ell$.
The group $\hat W$ acts transitively on $\Pc$, 
because $\Delta_{(w\ell)}=w(\Delta_{(\ell)})$ for all $w\in\hat W$.
The stabilizer of $(\ell_0)$ in $\hat W$ 
is $\hat W_{(\ell_0)}\hat W_{\ell_0}$.
If $s_{\hat\b}\in\hat W_{\ell_0}$ then $y_{\hat\b^\vee}(\ell_0)=1$,
hence $\b\in\Delta_{(\ell_0)}.$
Thus
$\hat W_{\ell_0}\subseteq\hat W_{(\ell_0)}$,
because
$\hat W_{\ell_0}$ is generated by reflections.
Therefore the stabilizer of $(\ell_0)$ in $\hat W$ 
is $\hat W_{(\ell_0)}$. 
Set ${}^{\{\ell_0\}}\Sb$
(resp. ${}^{\{\ell_0\}}\Sb_{(\ell)}$)
equal to the projective limit
$\pro\,\Sb/\langle E\rangle^n$
with $E\subset\hat W\ell_0$ 
(resp. $E\subset(\ell)$)
finite and $n\geq 0$. 
Hence ${}^{\{\ell_0\}}\Sb\simeq\prod_{(\ell)\in\Pc}{}^{\{\ell_0\}}\Sb_{(\ell)}$.
The tensor product 
${}^{\{\ell_0\}}\Hb={}^{\{\ell_0\}}\Sb\otimes_{\Sb}\Hb$ is a ring.
The ring ${}^{\{\ell_0\}}\Sb_{(\ell)}$ 
is a direct summand in ${}^{\{\ell_0\}}\Sb$.
The identity in ${}^{\{\ell_0\}}\Sb_{(\ell)}$ 
is identified with an idempotent in ${}^{\{\ell_0\}}\Sb$, 
denoted by $e_{(\ell)}$.
The same computations as in \cite{L, 8.13-16} 
yield a ${}^{\{\ell_0\}}\Sb_{(\ell_0)}$-linear ring isomorphism 
$${}^{\{\ell_0\}}\Hb_{(\ell_0)}\to 
e_{(\ell_0)}\cdot{}^{\{\ell_0\}}\Hb\cdot e_{(\ell_0)}\ \roman{ such\ that\ }
t_w\mapsto e_{(\ell_0)}t_w e_{(\ell_0)},\ \forall w\in\hat W_{(\ell_0)}.$$
By Proposition 7.2 and $(A)$ we have a chain of equivalences
$$
\M_\Pc\bigl({}^{\{\ell_0\}}\Hb_{(\ell_0)}\bigr)-mod^\infty
\to
{}^{\{\ell_0\}}\Hb_{(\ell_0)}{\ts -}\mod^\infty 
\to
{}^{\{\l_0\}}\Hb_{(\ell_0)}{\ts -}\mod^\infty.
$$
The rings 
${}^{\{\ell_0\}}\Hb$, ${}^{\{\ell_0\}}\Hb_{(\ell_0)}$, 
${}^{\{\l_0\}}\Hb_{(\ell_0)}$ 
are endowed with the topologies induced by the corresponding inverse systems,
and $mod^\infty$ is the category of smooth finitely generated modules,
see 7.2 for the other notations.
The restriction ${}^{\{\ell_0\}}\Hb$-$mod\to\Hb$-$mod$ yields an equivalence 
${}^{\{\ell_0\}}\Hb$-$mod^\infty\to{}^{\{\ell_0\}}\Oc$.
Similarly,
the restriction ${}^{\{\l_0\}}\Hb_{(\ell_0)}$-$mod\to\Hb'_{(\ell_0)}$-$mod$ 
yields
an equivalence
${}^{\{\l_0\}}\Hb_{(\ell_0)}$-$\mod^\infty\to {}^{\{\l_0\}}\Oc_{(\ell_0)}$.
Thus it suffices to prove that the categories
$\M_\Pc({}^{\{\ell_0\}}\Hb_{(\ell_0)})$-$mod^\infty$
and
${}^{\{\ell_0\}}\Hb$-$mod^\infty$
are equivalent.

For each class $(\ell)\in\Pc$ 
we fix an element $w_{(\ell)}\in\hat W$ such that
$(\ell)=w_{(\ell)}(\ell_0)$ and 
$e_{(\ell)}\varphi_{w_{(\ell)}}\in{}^{\{\ell_0\}}\Hb\cdot e_{(\ell_0)}$,
see \cite{L, 8.16}.
We write $\varphi_{(\ell)}$ for $\varphi_{w_{(\ell)}}$.
Let $E_{(\ell)(\ell')}(h)\in\M_\Pc({}^{\{\ell_0\}}\Hb_{(\ell_0)})$
be the matrix with $(\ell),(\ell')$-th entry equal to $h$ and all 
other entries equal to zero.
The linear map
$${}^{\{\ell_0\}}\Hb\to\M_\Pc\bigl({}^{\{\ell_0\}}\Hb_{(\ell_0)}\bigr),\quad
\varphi_{(\ell)}h\varphi_{(\ell')}^{-1}\mapsto
E_{(\ell)(\ell')}(h),$$
is an embedding of topological rings with a dense image, see \cite{L, 8.16}.
The restriction yields the desired equivalence
$\M_\Pc({}^{\{\ell_0\}}\Hb_{(\ell_0)})$-$\mod^\infty\to
{}^{\{\ell_0\}}\Hb$-$mod^\infty$.

\vskip2mm

The last claim in $(iii)$ is proved as in $(B)$, 
see also \cite{L, Section 8}.
\qed

\subhead 2.6. Intertwiners in $\Ocu$\endsubhead
Assume that $\k$ is a field.
For any $J\subsetneq\hat I$ let
$\Hb_J=\Oplus_{w\in W_J}t_w\Sb\subset\Hb$.
It is a subring. 
We write $\Hbu$ for $\Hb_I,$ $\Ocu$ for $\Hbu$-$\mof$,
and $\underline{[\ell]}$ (or $\underline{[\ell]}_\k$ if necessary)
for $[\ell]_{W_\ell,\k}$.
Set $\Sb_\ell=\Sb/\underline{[\ell]}$,
$\Pu(\ell)=\Hbu\otimes_\Sb\Sb_\ell$,
and $1_\ell=1\otimes 1\in\Pu(\ell)$.

Let ${}^\ell\Hbu$ be 
the specialization of $\Hbu$ at the central character 
$W\ell\in\Spec(\Sb^W)$. 
Set ${}^\ell\Ocu={}^\ell\Hbu$-$\mof$.
The module $\Pu(m)$ lies in ${}^\ell\Ocu$ for all $m\in W\ell$ because
$\la\ell\ra^W\subseteq\underline{[m]}$ and
$\la\ell\ra^W$ lies in the center of $\Hbu$.
It is projective (as in Proposition 2.2$(i)$).

For each $w\in W$ the intertwining operator 
$\Phi_w(\ell)\,:\,\Pu(w\ell)\to\Pu(\ell)$ is the unique 
$\Hbu$-homomorphism taking
$1_{w\ell}$ to $\phi_w 1_\ell$.
The same argument as for Lemma 2.3 implies that 
$\Phi_{s_i}(\ell)$ is invertible if and only if
$y_{\a_i^\vee}(\ell)\neq\zeta_i^{\pm 1}.$

The connected components of the set 
$X_\RR\setminus\bigcup_{\b^\vee\in\Delta^\vee}H_{\b^\vee}$ are the chambers. 
Let $C_\pm$ be the chamber containing $\pm\rho$, and
$C_w=\{w^{-1}\mu\,;\,\mu\in C_+\}$.

Set
$\Hc_\ell=\{\b^\vee\in\Delta^\vee\,;
\,y_{\b^\vee}(\ell)=\zeta^{\pm 1}_{\b^\vee}\}$,
and 
$U_\ell=X_\RR\setminus\bigcup_{\b^\vee\in\Hc_\ell}H_{\b^\vee}$.
The group $W_\ell$ acts on $U_\ell$.
A domain is a minimal subset in $U_\ell$
containing a connected component and stable by $W_\ell$.
Let $\Du_w$ be the unique domain containing $C_w$, 
and $\Dcu$ be the set of domains. 

\proclaim{Proposition}
(i)
$\Pu(w_1\ell)$, $\Pu(w_2\ell)$ are isomorphic 
whenever $\Du_{w_1}=\Du_{w_2}$.

(ii)
Assume that $\hat W_\l=W_\ell$.
There is a unique injection $\dagger\,:\,\Dcu\to\Dc$ such that 
$\Du_{w_2}^\dagger=D_{w_1}$ if $w_2v=x_\kappa ww_1$ 
and $w\in W$, $v\in\hat W_\l$, $\kappa\in Y$ far enough inside $C_+$. 
\endproclaim

\noindent{\sl Proof.}
Claim $(i)$ is immediate using the condition for
the invertibility of the intertwining operator given above.
Claim $(ii)$ is easy and is left to the reader.
\qed

\head 3. Reminder on Knizhnik-Zamolodchikov trigonometric connection\endhead

This section contains standard results 
on Knizhnik-Zamolodchikov trigonometric connection.
See \cite{GGOR} for the analogue in the rational case.
The sheaf of regular functions on a scheme $S$ is denoted by $\Oc_S$.
If $S$ is smooth, the sheaf of differential operators on $S$ is denoted 
by $\Dc_S$.

\subhead 3.1\endsubhead
Assume that $\k$ is a field.
Set $T_\reg=\{x_\b\neq 1\,;\,\forall\b\in\Delta\}\subset T$.
Let $\Db_\reg$ be the ring of algebraic differential operators on $T_\reg$.
For each $j\in I$ set 
$$D_j=\partial_{\xi_j}-
\sum_{\b\in\Delta_+}h_\b\b_j\vartheta_\b+\tilde\rho_j\in\Db_\reg
\ \roman{with}\ 
\tilde\rho={1\over 2}\sum_{\b\in\Delta_+}h_\b\otimes\b.
\leqno(3.1.1)$$
Put $\Rb_\reg=\k[T_\reg]$, 
and $\Hb'_\reg=\Rb_\reg\otimes_{\Rb}\Hb'$.
Set $\theta_\b=(1-x_{-\b})^{-1}\otimes(1-s_\b)\in\Hb'_\reg$.

\proclaim{Lemma}
(i)
There is a unique $\k$-algebra structure on $\Hb'_\reg$ extending $\Hb'$. 

(ii)
There is a unique ring isomorphism
$\Db_\reg\rtimes\k W\to\Hb'_\reg$
such that 
$\partial_{\xi_j}\mapsto\nabla_j
:=\xi_j+\sum_{\b\in\Delta_+}h_\b\b_j\theta_\b-\tilde\rho_j$,
$f\mapsto f$, $w\mapsto w$
for all $j\in I$, $w\in W$ and $f\in\Rb_\reg$.
\endproclaim

\noindent{\sl Proof.}
$(i)$ is immediate because $\Hb'$ is a free left $\Rb$-module and 
$\{f\in\Rb\,;\,f(T_\reg)=0\}$ is a left Ore set in $\Hb'$ by (2.1.2).

Let us check $(ii)$.
Observe that $D_j$ preserves the subspace $\Rb\subset\Rb_\reg$.
Identifying $\Rb$ with the module $\Hb'\otimes_{\Hbu'}\k$
induced from the trivial representation of $\Hbu'$ on $\k$, we get
a representation of $\Hb'$ on $\Rb$ such that
$w(g)={}^w g$, $\xi_j(g)=D_j(g)$ and $f(g)=fg$ for each
$f,g\in\Rb$ and $w\in W$.
This action extends obviously to an action of $\Hb'_\reg$ on $\Rb_\reg$.
Hence there is a ring homomorphism $\Hb'_\reg\to\Db_\reg\rtimes\k W$
such that $\xi_j\mapsto D_j$, $f\mapsto f$, $w\mapsto w$.
It is obviously surjective.
It is also injective because
the representation of $\Hb'$ on $\Rb$ above is faithful, 
by a well-known lemma of Cherednik.
\qed

\subhead 3.2\endsubhead
For each $\Hb'$-module $M$ we set $M_\reg=\Rb_\reg\otimes_{\Rb}M$.
Composing the localization $\Oc'\to\Hb'_\reg$-$\Modc$, $M\mapsto M_\reg$,
the isomorphism 3.1,
and the sheafification 
$\Db_\reg$-$\Modc\to\Dc_{T_\reg}$-$\Modc$,
we get a functor
$\Lc\,:\,\Oc'\to\Dc_{T_\reg}\rtimes\k W$-$\Modc$.
For any $M$ in $\Oc'$ the 
$\Dc_{T_\reg}\rtimes\k W$-module $\Lc(M)$ 
is locally free of finite rank over $\Oc_{T_\reg},$ 
because $\Lc(M)$ is a $\Dc_{T_\reg}$-module which is
coherent over $\Oc_{T_\reg}$ 
(since $M$ is finitely generated over $\Rb$).

\subhead 3.3\endsubhead
Set $\k=\CC$.
Let $z_i,$ $i\in I$, be the obvious coordinates on $\CC^I$.
For any $\b\in\Delta$ we write $z^\b$ for $\prod_iz_i^{\b_i}$.
Let $D_\infty\subset\CC^I$ be the divisor $\{\prod_{i\in I}z_i=0\}$.
The map $(x_{\a_i})\,:\,T\to\CC^I$ is an isomorphism
onto $\CC^I\setminus D_\infty$. 
Set $D_\Delta=\bigcup_{\b\in\Delta}\{z^\b=1\}$,
and $D=D_\infty\cup D_\Delta$.
Then $T_\reg$ is identified with the open set
$\CC^I_\reg=\CC^I\setminus D$. 

Let $\CC^I\to\CC^I_\reg/W$, $u\mapsto[u]$ be the obvious projection.
Fix $\copyright\in (0,1)^I$, 
and $\l^\vee_c\in X^\vee_\CC$ such that $\e^{\l_c^\vee}=\copyright$.
The fundamental group $\Pi_1(\CC^I_\reg/W,[\copyright])$ 
is generated by the homotopy classes of the paths
$\g_j,\tau_j\,:\,[0,1]\to\CC^I_\reg/W$ such that
$$\g_j(t)=[\copyright\cdot\e^{t\o^\vee_j}],\quad
\tau_j(t)=[\copyright\cdot\e^{-t(\a_j\,:\,\l_c^\vee)\a_j^\vee}].$$
It is isomorphic to the affine braid group $B_{\hat W}$ associated to $\hat W$, 
see \cite{H, \S 2} for more details and references.

From now on we assume that $\k=\A$, $\F$ or $\CC$.
For any finite dimensional $\CC$-vector space $V$ 
we call holomorphic function $\CC^I\to V((\varpi))$ a formal series
$\sum_{n\gg -\infty}a_n\varpi^n$ where each $a_n$
is a holomorphic function $\CC^I\to V.$

Given a $W$-equivariant $\k$-vector bundle $V$ over $\CC^I_\reg$
with a $W$-invariant integrable connection $\nabla$, 
let $V^\nabla$
be the set of $W$-invariant holomorphic horizontal sections
of $V$ over the simply connected cover 
$\tilde\CC^I$ of $\CC^I_\reg$.
It is a free $\k$-module of rank equal to the rank of $V$.

The group $B_{\hat W}$ acts on $V^\nabla$ by monodromy.
The functor $V\mapsto V^\nabla$ is exact,
from the category of $W$-equivariant vector bundles on $\CC^I_\reg$
with a $W$-invariant integrable connection to $\k B_{\hat W}$-$\mof$.
It restricts to an equivalence from 
the category of $W$-equivariant vector bundles on $\CC^I_\reg$ 
with a regular integrable $W$-invariant connection to $\k B_{\hat W}$-$\mof$.

If $\k=\A$ we have $\CC\otimes_\A V^\nabla=(\CC\otimes_\A V)^\nabla$ and
$\F\otimes_\A V^\nabla=(\F\otimes_\A V)^\nabla$.

\subhead 3.4\endsubhead
Let $\Mc\,:\,\Oc'\to\k B_{\hat W}$-$\mof$ 
be the functor $M\mapsto\Lc(M)^\nabla$.

\proclaim{Lemma}
Fix $M,N\in\Oc'$.

(i)
The canonical map 
$\Hom_{\Oc'}(M,N)\to\Hom_{B_{\hat W}}(\Mc(M),\Mc(N))$ is injective
if $N$ is torsion-free over $\Rb$.

(ii)
The $\Dc_{\CC^I_\reg}$-module $\Lc(M)$ 
has regular singularities along $D$.
\endproclaim

\noindent{\sl Proof.}
The restriction 
$\Hom_{\Hb'}(M,N)\to\Hom_{\Hb'_\reg}(M_\reg,N_\reg)$
is injective if $N$ is torsion free over $\Rb$.
Assigning to a horizontal section on $\CC^I_\reg/W$
its value in the fiber at a given point is an injective map.
Thus the map
$$\Hom_{\Oc'}(M,N)\to\Hom_{B_{\hat W}}(\Mc(M),\Mc(N))$$
is injective.
Claim $(i)$ is proved.

Fix a $\Hbu'$-module $M$ in $\Ocu'$.
The horizontal sections of $\Lc\Ic(M)$ are 
the elements in $\Rb_\reg\otimes_\k M$ 
annihilated by the operator $\nabla_j$ for all $j\in I$,
see Lemma 3.1.
Using (2.1.2) we get
$$\nabla_j=
\partial_{\xi_j}\otimes 1+1\otimes\xi_j
-\sum_{\b\in\Delta_+}h_\b\b_j\theta_\b-\tilde\rho_j.$$
Hence the elements of $\Mc\Ic(M)$ are the $W$-invariant 
maps $\CC^I_\reg\to M$ which are annihilated by the connection 
$\d-\sum_jA_j\d z_j/z_j$, with 
$$A_j=\tilde\rho_j-\xi_j
-\sum_{\b\in\Delta_+}h_\b{\b_j\,z^\b\over 1-z^\b}(1-s_\b).
$$
This connection is the trigonometric Knizhnik-Zamolodchikov connection 
on the vector bundle $\CC^I_\reg\times_W M_\reg$ over $\CC^I_\reg/W$. 
It has regular singularities along $D$ and at infinity.

The category of $\Oc_{\CC^I_\reg}$-coherent $\Dc_{\CC^I_\reg}$-modules
with regular singularities is stable by subquotients.
Therefore $\Lc(M)$ has regular singularities for each $M\in{}^\l\Oc$
by Proposition 2.2.$(ii)$.

The category of $\Oc_{\CC^I_\reg}$-coherent $\Dc_{\CC^I_\reg}$-modules
with regular singularities is stable by extensions. 
Therefore $\Lc(M)$ has regular singularities for each $M\in{}^{\{\l\}}\Oc$.
Then, $(ii)$ follows from Proposition 2.1.$(ii)$.
\qed

\vskip3mm

\noindent{\bf Notation.}
If $M\in\Ocu'$ we write $M^\nabla$ for $\Mc\Ic(M)$. 

\head 4. Monodromy\endhead

We fix a branch of the logarithm.
Put $z^a=\exp(a\log(z))$ for any $a$.
Set $\k=\CC$. 
Fix $\l_0\in X_\CC$ such that $\hat W_{\l_0}\subseteq W$.
Set $\ell_0=\e^{\l_0}$ and $\zeta^{1/2}_{0i}=\e^{h_{0i}/2}$.
Note that $\hat W_{\l_0}=W_{\ell_0}$ by Lemma 1.3, 
hence $\hat W_{\l_0}$ is generated by reflections.
We assume that $\zeta_{0i}\neq 1,-1$ for each $i$.

\subhead 4.1\endsubhead
The modules $\Pu(\mu_0)^\nabla$, with $\mu_0\in X_\CC$,
have been studied by several authors
when the parameters are generic, 
see \cite{C3, Proposition 3.4} for instance.
It is important, for us, to have precise information 
for non generic values of the parameters.

\proclaim{Theorem}
(i)
$\Pu(\hat w\l_0)^\nabla=\Pu(w\ell_0)$ 
for all $\hat w\in\hat W$, $w\in W$ 
with $D_{\hat w}=\Du_{w}^\dagger$.

(ii)
$\Mc$ factors through a functor $\Oc'\to\Ocu$.

(iii)
$\Mc$ is fully faithful on $\Ic(\Ocu')$.
\endproclaim

\noindent{\sl Proof of $(i)$.}
Fix $\mu_0\in\hat W\l_0$ and $m_0=\e^{\mu_0}$. 
The computation of $\Pu(\mu_0)^\nabla_\CC$ 
uses a reduction to the rank one case as in \cite{C3}.
To do so, we first deform $\Pu(\mu_0)_\CC$ over $\A$.
Then we fix a fundamental matrix solution over $\F$. 
From now on $\k=\A,\F$ or $\CC$.

\vskip2mm

\noindent $(A)$
Set
$$X_0=\{\epsilon\in X_\CC\,;\,
(w\epsilon)_j\neq (w'\epsilon)_j,\,\forall w\neq w'\in W,\,\forall j\in I\}.$$
Put $\mu=\mu_0+\varpi\epsilon$, with $\epsilon\in X_0$.
Set $Q=\hat W_{\mu_0}\cdot\mu.$
From now on let $\nu$ denote any element in $Q$.
Set $\Sb_{Q,\A}$ as in 2.1. 
The ring $\Sb_{Q,\A}$ is local.
Let $\mb_{Q,\A}$ be the maximal ideal. 
Set $\Sb_{Q,\k}=\k\otimes_\A\Sb_{Q,\A}$ 
for $\k=\F$ or $\CC$.
We claim that $\Sb_{Q,\CC}=\Sb_{\mu_0,\CC}$.
We have $(\Sb'_\A/[\mu]_{\hat W_{\mu_0},\A})\otimes_\A\CC=\Sb_{\mu_0,\CC}$ 
by Remark 2.1.
Hence the obvious surjective map 
$\Sb'_\A/[\mu]_{\hat W_{\mu_0},\A}\to\Sb_{Q,\A}$
specializes to a surjective map $\Sb_{\mu_0,\CC}\to\Sb_{Q,\CC}$.
The claim follows, because
$\dim(\Sb_{\mu_0,\CC})=\hat n_{\mu_0}$ by Chevalley's theorem, 
and $\dim(\Sb_{Q,\CC})=\hat n_{\mu_0}$ because 
$\dim(\Sb_{Q,\F})=\hat n_{\mu_0}$, since $\hat W_\nu=\{1\}$,
and $\Sb_{Q,\A}$ is free over $\A$ because
$\Sb_{Q,\A}\subset\Oplus_\nu\Sb_{\nu,\A}$ and 
$\Sb_{\nu,\A}=\A$.

Put $\Pu=\Hbu'\otimes_{\Sb'}\Sb_Q$. 
The module $\Pu$ lies in $\Ocu'$ and $\Pu_\CC=\Pu(\mu_0)_\CC$.
Let $Y_j,$ $T_j$ be the monodromy operators on $\Pu^\nabla$
along $\g_j$, $\tau_j$ respectively.
The assignement 
$y_j\mapsto\e^{\tilde\rho_j} Y_j,$
$t_j\mapsto\zeta^{1/2}_{0j} T_j$ 
extends uniquely to a representation of $\Hbu_\F$ on $\Pu_\F^\nabla$
by \cite{C1, Proposition 8}.
The canonical maps
$\F\otimes_\A\Pu^\nabla_\A\to\Pu^\nabla_\F$ and 
$\CC\otimes_\A\Pu^\nabla_\A\to\Pu^\nabla_\CC$
commute to the $B_{\hat W}$-action.
Therefore the representation of $B_{\hat W}$ on $\Pu^\nabla_\k$ 
factors also through $\Hbu_\k$ if $\k=\A,\CC$.

\vskip2mm

\noindent $(B)$
Assume that 
$$(\mu_0\,:\,\b^\vee)\in\RR_{\ll0}+i\RR,
\qquad\forall\b^\vee\in\Delta^\vee_+.\leqno(4.1.1)$$
We first prove that $\Pu^\nabla_\CC$ is cyclic over $\Hbu_\CC$.
Then we prove that $\Pu^\nabla_\CC\simeq\Pu(m_0)_\CC$.

Set $\psi_w=\phi'_w\otimes 1\in\Pu$ for each $w\in W$.
Hence $\psi_w\in w\psi_1\pi_w+\sum_{w'<w}w'\psi_1\Sb_Q$,
where $\pi_w$ is the product of all $\xi_{\a^\vee}$ with 
$\a^\vee\in\Delta_+^\vee\cap w^{-1}\Delta_-^\vee$.
The image of $\pi_w$ in $\Sb_{Q,\A}$ is invertible :
we have $\pi_w\notin\mb_{Q,\A}$ 
because the image of $\pi_w$ in $\Sb_{\mu_0,\CC}$ does not lie
in $\la\mu_0\ra_\CC$ since $\mu_0$ is regular by (4.1.1).
Thus $(\psi_w)$ is a $\Sb_{Q,\A}$-basis of $\Pu_\A$.

The obvious right $\Sb_Q$-action on $\Pu$ 
commutes to the left $\Hbu'$-action,
thus $A_j$ (= the connection matrix in 3.4) is $\Sb_Q$-linear, 
hence $\Pu^\nabla$ is a $(\Hbu,\Sb_Q)$-bimodule.
If $k\in\ZZ$ is non-zero, the image of the element
$t=k+{}^{w^{-1}}\!\!\xi_j-{}^{{w'}^{-1}}\!\!\xi_j$ in 
$\Sb_{Q,\F}$ is invertible,
because $\Sb_{Q,\F}=\Oplus_\nu\Sb_{\nu,\F}$ 
and the projection of $t$ in $\Sb_{\nu,\F}$ is invertible 
(since ${}^{w^{-1}}\!\!\xi_j\in(w\nu)_j+\la\nu\ra$ 
and $\epsilon\in X_0$).
Put $A_{j0}=\tilde\rho_j-\xi_j$.
The element $A_{j0}\in\Sb'$ 
is identified with its projection in $\Sb_Q$ whenever needed.
Set $z^{A_0}=\Prod_jz_j^{A_{j0}}$.
There is a unique $\Sb_{Q,\F}$-basis $(\psi_w^\nabla)$ of
$\Pu^\nabla_\F$ such that the function
$$(z_j)\mapsto\psi^\nabla_w(z_j)\cdot z^{-B},$$
where $B={}^{w^{-1}}\!\!\!A_0,$
is holomorphic on $\CC^I\setminus D_\Delta$ and equals $\psi_w$ at 0.
By Proposition 7.1 there is also a $\Sb_{Q,\A}$-basis 
$(b^\nabla_w)$ of $\Pu_\A^\nabla$ such that 
$$b^\nabla_w\in\psi_w^\nabla+
\sum_{w'\mu_0<w\mu_0}\psi_{w'}^\nabla\cdot\Sb_{Q,\F}.\leqno(4.1.2)$$


Let $b_w$ be the image of $b_w^\nabla$
by the unique $\Sb_{Q,\F}$-linear isomorphism
$\Pu^\nabla_\F\to\Pu_\F$ 
such that $\psi_w^\nabla\mapsto\psi_w$.
Let $\Pu^\nabla_\A$ denote also the $\Sb_{Q,\A}$-span of $(b_w)$. 
The monodromy $\Hbu_\F$-action on $\Pu_\F^\nabla$
yields a representation of $\Hbu_\A$ on $\Pu_\F$ 
which preserves $\Pu^\nabla_\A$. 
For each $\eta_0,\eta'_0\in X_\CC$ we write 
$\eta_0\succ\eta'_0$ if 
$\eta_0-\eta'_0\in Y_{\RR,+}\setminus\{0\}+iX_\RR$.
Note that $w\mu_0\succ w'\mu_0$ if $w>w'$, by (4.1.1),
or if $w\mu_0>w'\mu_0$.
We claim that
$$t_w b_1\in b_w\cdot\Sb_{Q,\A}^\times+
\sum_{w'\mu_0\prec w\mu_0}b_{w'}\cdot\Sb_{Q,\A},\ \roman{and}\ 
\Sb_\A b_1=b_1\cdot\Sb_{Q,\A}.\leqno(4.1.3)$$
The proof is given in $(D)$.
Then an easy induction implies that
$\Pu_\A^\nabla=\Hbu_\A b_1$, 
hence that $\Pu_\CC^\nabla=\Hbu_\CC b_1$.

The series $\e^{\xi_j}=\sum_{k\geq 0}(2i\pi\xi_j)^k/{k!}$
converges in $\Sb_{Q,\A}$,
because $\xi_j\in(\mu_0)_j+\mb_{Q,\A}$ and $\mb_{Q,\A}$
is pronilpotent.
Let $\Sb_\A\to\Sb_{Q,\A}$, $p\mapsto p(\e^\xi)$ 
be the ring homomorphism such that $y_j\mapsto\e^{\xi_j}$.
It is surjective. 
The map $p\mapsto p(\e^\xi)$ 
factors through a ring isomorphism $\Sb_{m_0,\CC}\to\Sb_{\mu_0,\CC}$
(fix $\kappa\in Y$ such that $\mu_0+\kappa\in W\l_0$;
then $\hat W_{\mu_0}=x_\kappa^{-1}W_{m_0}x_\kappa,$
because $\hat W_{\mu_0}=x_\kappa^{-1}\hat W_{\mu_0+\kappa}x_\kappa$
and $\hat W_{\l_0}=W_{\ell_0}$;
the claim follows since $\dim(\Sb_{m_0,\CC})=\dim(\Sb_{\mu_0,\CC})$). 
In particular $\underline{[m_0]}b_1=0$ by (4.1.4), because $[\mu_0]b_1=0$.
Therefore there is a unique surjective $\Hbu_\CC$-linear map
$\Pu(m_0)_\CC\to\Pu^\nabla_\CC$ such that $1_{m_0}\mapsto b_1$.
It is invertible because both modules 
have the same dimension over $\CC$.

\vskip2mm

\noindent $(C)$
Fix $\hat w\in\hat W$, $w\in W$ as in $(i)$.
We may assume that $wv=x_\kappa w'\hat w$ 
with $w'\in W$, $v\in\hat W_{\l_0}$, and $\kappa\in Y$ far inside $C_+$, 
because $D_{\hat w}=\Du_{w}^\dagger$.
In particular the alcove $A_{\hat w}$ 
is far inside $D_{\hat w}$.
Put $\mu_0=w'\hat w\l_0$.
Then (4.1.1) holds.
Thus $\Pu^\nabla_\CC=\Pu(m_0)_\CC$ by $(A)$.
We have also $m_0=w\ell_0$ because $\mu_0+\kappa=w\l_0$.
Thus $\Pu(\hat w\l_0)^\nabla_\CC=\Pu(w\ell_0)_\CC$, because
$\Phi'_{w'}(\hat w\l_0)\,:\,P(\mu_0)_\CC\to P(\hat w\l_0)_\CC$
is invertible (since $A_{\hat w}$ is far inside $D_{\hat w}$).

\vskip2mm

\noindent $(D)$
Let $G\,:\,\CC^I_\reg\to\End(\Pu_\F)$
be the fundamental matrix solution such that $G\psi_w=\psi_w^\nabla$.
We have $G=Hz^{A_0}$ with 
$H\,:\,\CC^I\setminus D_\Delta\to\End(\Pu_\F)$
holomorphic such that $H(0)=\Id$, and
$Y_j=G(z\e^{\o_j^\vee})^{-1}G(z),$
$T_j=G(s_jz)^{-1}s_jG(z).$
Thus 
$${}^wp\,\psi=\psi\cdot p(\e^\xi),\quad\forall p\in\Sb_\F,\,\forall w.
\leqno(4.1.4)$$
Hence the second part of (4.1.3) is immediate.

Let us prove the first part.
We first claim that for each $w\in W$ 
there is an invertible element $p_w\in\Sb_{Q,\F}$ such that
$$t_w\psi_1\in\psi_w\cdot p_w+\sum_{w'<w}\psi_{w'}\cdot\Sb_{Q,\F}.
\leqno(4.1.5)$$
To do so, observe that
$$\psi_w\cdot\Sb_{Q,\F}=\{\psi\in\Pu_\F\,;\,
{}^wp\,\psi=\psi\cdot p(\e^\xi),\,\forall p\in\Sb_\F\}.\leqno(4.1.6)$$
Indeed, the direct inclusion is immediate, while
the inverse one holds because $\Pu_\F=\Oplus_{w'}\psi_{w'}\cdot\Sb_{Q,\F}$
and, for each $w\neq w'$, there is an element $p\in\Sb_\F$ such that
$${}^wp(\e^\xi)-{}^{w'}p(\e^\xi)\in\Sb_{Q,\F}^\times$$
(because
$\Sb_{Q,\F}=\Oplus_\nu\Sb_{\nu,\F}$, 
and there is $p$
such that $p(w\e^{\nu})\neq p(w'\e^{\nu})$ 
for each $\nu$ 
since $W_{\e^\nu}=\{1\}$).
Then, (4.1.6) implies that $\phi_w\psi_1\in\psi_w\cdot\Sb_{Q,\F}$, 
and (4.1.5) follows.

Using (4.1.2) and (4.1.5) we get
$$t_wb_1\in b_w\cdot p_w+\sum_{w'\mu_0\prec w\mu_0}b_{w'}\cdot\Sb_{Q,\A}
\leqno(4.1.7)$$
for some $p_w\in\Sb_{Q,\A}\cap\Sb_{Q,\F}^\times$.
We must prove that $p_w\in\Sb_{Q,\A}^\times$.
We prove it by induction on $\ell(w)$.
Fix $v\in W$ such that $s_jv>v$.
By (4.1.5) there is an element $q\in\Sb_{Q,\F}$ such that
$$t_j\psi_v\in\psi_{s_jv}\cdot q+
\sum_{v'<s_jv}\psi_{v'}\cdot\Sb_{Q,\F}.
\leqno(4.1.8)$$
By (4.1.2) and (4.1.7) we have $q\in\Sb_{Q,\A}$.
It is sufficient to prove that
$q\in\Sb_{Q,\A}^\times$.
To simplify we write $j$ for $\{j\}$ 
and $P_j$ for $\Hb'_j\otimes_{\Sb'}\Sb_{vQ}$,
where $\Sb_{vQ}$ is defined as $\Sb_Q$ in $(A)$.
From now on $w$ is either $v$ or $s_jv$.
Set $\varphi_w=\phi'_{wv^{-1}}\otimes 1\in P_j$.
Then $(\varphi_w)$ is a $\Sb_{vQ}$-basis of $P_j$.

Let $P_j^\nabla$ be the set of holomorphic functions 
$f\,:\,\CC\setminus\{0,1\}\to P_j$ 
such that 
$$z_j\partial_{z_j}f-A_{j0}f+h_{0j}z_j{1-s_j\over 1-z_j}f=0.$$
It is a right $\Sb_{vQ}$-module.
Let $Y'_j$, $T'_j$ the monodromy operators around 0 and 1.
Since $y_k$ lies in the center of $\Hb_j$
if $k\neq j$, the assignement
$y_j\mapsto\e^{\tilde\rho_j}Y'_j$,
$t_j\mapsto\zeta^{1/2}_{0j}T'_j$ 
extends to a representation of $\Hb_j$ on $P_j^\nabla$ such that
$y_k m=m\cdot\e^{{}^{v^{-1}}\!\!\xi_k}$ 
for each $k\neq j$ and each $m\in P_j^\nabla$. 
Let $G_j$ be the fundamental matrix solution such that $G_j=H_jz_j^{A_{j0}}$ 
with $H_j\,:\,\CC\setminus\{1\}\to\End(P_{j,\F})$ holomorphic and $H_j(0)=\Id$. 
Set $\varphi_w^\nabla=G_j\varphi_w.$
There is a unique $\Sb_{vQ,\F}$-linear isomorphism 
$P_{j,\F}^\nabla\to P_{j,\F}$ such that 
$\varphi_w^\nabla\mapsto\varphi_w.$
It yields a representation of $\Hb_{j,\F}$ on $P_{j,\F}$.

Let $\theta_j\,:\,P_j\to\Pu$ be the $\Hb'_j$-linear map 
such that $\varphi_w\mapsto\psi_w$.
Note that $\theta_j(m\cdot{}^vp)=\theta_j(m)\cdot p$ for each $m\in P_j$,
$p\in\Sb_Q$.
We have 
$$\theta_j(t_j\varphi_v)=
\lim_{\eps\to 0}\eps^{D_j}\circ t_j\circ\eps^{-D_j}(\psi_v)
\ \roman{with}\ D_j=\sum_{k\neq j}A_{k0},$$ 
because
$\theta_j\circ G_j=
\lim_{\eps\to 0}(G\eps^{-D_j})|_{C_\eps}\circ\theta_j$
with
$C_\eps=\bigcap_{k\neq j}\{z_k=\eps\}\subset\CC^I$.
Thus $t_j\varphi_v=\varphi_{s_jv}\cdot{}^vq$ modulo 
$\varphi_v\cdot\Sb_{vQ}$, with $q$ as in (4.1.8), because
$D_j(\psi_w)=\psi_w\cdot a$ for some element $a\in\Sb_Q$ 
which is independent on $w\in\{v,s_jv\}$,
and because $t_j\psi_v$ is a linear combination of the elements
$\psi_{v'}$ with $v'\leq s_jv$ by (4.1.8).
Therefore to prove (4.1.3) it suffices to check that 
$$t_j\varphi_v\in\varphi_{s_jv}\cdot\Sb_{vQ,\A}^\times+
\varphi_v\cdot\Sb_{vQ,\F}.$$

Since $\Sb_{vQ,\A}\subseteq\Oplus_\nu\Sb_{v\nu,\A},$
an element in $\Sb_{vQ,\A}$ is invertible 
if and only if its image in $\Sb_{v\nu,\A}$ is invertible.
There is a unique $\Hb_j'$-linear map $P_j\to P_j(v\nu)$ 
taking $1\otimes 1$ to $1\otimes 1$.
It commutes to the right actions of 
$\Sb_{vQ}$ on $P_j$, and of $\Sb_{v\nu}$ on $P_j(v\nu)$.
Let $\bar\varphi_w$ be the image of $\varphi_w$.
Since $\Sb_{v\nu,\A}=\A$ for each $\nu$,
it is enough to prove that
$$t_j\bar\varphi_v\in\bar\varphi_{s_jv}\cdot\A^\times
+\bar\varphi_v\cdot\F.\leqno(4.1.9)$$
Let $\Gamma$ be the gamma function.
For each $z\in\CC+\varpi\A^\times$ set 
$a(z)=(\zeta_{0j}^{1/2}-\zeta_{0j}^{-1/2})(\e^z-1)^{-1}$ and
$b(z)=\Gamma(z)\Gamma(1+z)\Gamma(h_{0j}+z)^{-1}\Gamma(1-h_{0j}+z)^{-1}$.
Then \cite{C1, Theorem 10} yields 
$$t_j\bar\varphi_v=\bar\varphi_{s_jv}\cdot b(-\g)+\bar\varphi_v\cdot a(-\g),$$
with $\g=(v\nu\,:\,\a_j^\vee)$.
Note that $\g=(v\mu_0\,:\,\a_j^\vee)+\varpi(v\epsilon\,:\,\a_j^\vee)$, 
where $(v\mu_0\,:\,\a_j^\vee)\notin\{0,\pm h_{0j}\}+\ZZ_{\geq 0}$ 
by (4.1.1) because $s_jv>v$. 
Thus $b(-\g)\in\A^\times$, because
$\Gamma$ does not vanish anywhere,
and has a simple pole at each non positive integer.
Hence (4.1.9) holds.
\qed

\vskip3mm

\noindent{\sl Proof of $(ii)$. }
Set $\k=F$. 
Fix $\l\in X_\A$, $h_i\in\A$ such that 
$(\l,h_i)=(\l_0,h_{0i})$ modulo $\varpi$,
and $(\l,h_i)$ is generic over $\F$. 
Then $P(\mu)_\F=P(\l)_\F$ and
$\Pu(\l)^\nabla_\F=\Pu(\e^\l)_\F$ for all $\mu\in\hat W\l$.
Thus the $\F B_{\hat W}$-action on $\Mc(M)$ 
factors through $\Hbu$ for all $M$ in ${}^\l\Oc_\F$ by Proposition 2.2, 
yielding a functor $\Mc\,:\,{}^\l\Oc_\F\to{}^\ell\Ocu_\F$.

Fix $\k=\A$, and $(\l,h)$ as above.
For each $M$ in ${}^\l\Oc_\A$, 
$\Mc(M)$ is free over $\A$,
$\Mc(\F\otimes_\A M)=\F\otimes_\A\Mc(M)$, and 
$\Mc(\F\otimes_\A M)\in {}^\ell\Ocu_\F$.
Hence $\Mc(M)\in {}^\ell\Ocu_\A$,
thus $\Mc(\CC\otimes_\A M)=\CC\otimes_\A\Mc(M)\in {}^{\ell_0}\Ocu_\CC$.

Fix $\k=\CC$.
Then 
$\Mc({}^{\l_0}\Oc)\subset{}^{\ell_0}\Ocu$.
Therefore
$\Mc({}^{\{\l_0\}}\Oc)\subset{}^{\{\ell_0\}}\Ocu$,
because an object in ${}^{\{\l_0\}}\Oc$ 
has a filtration whose associated graded lies in ${}^{\l_0}\Oc$
and $\Mc$ is exact.
\qed

\vskip3mm

\noindent{\sl Proof of $(iii)$.}
Fix $M,N\in\Ocu'$.
Since $\Ic(N)$ is torsion-free over $\Rb$ the natural map
$$\Hom_{\Oc'}(\Ic(M),\Ic(N))\to
\Hom_{B_{\hat W}}(M^\nabla,N^\nabla)$$
is injective by Lemma 3.4.$(i)$.
The functor of horizontal sections yields an isomorphism
$$\Hom_{\Hb'_\reg}(\Ic(M)_\reg,\Ic(N)_\reg)\to
\Hom_{B_{\hat W}}(M^\nabla,N^\nabla)$$
by Lemma 3.1.$(ii)$, Lemma 3.4.$(ii)$.
We must check that the restriction map
$$\Hom_{\Hb'}(\Ic(M),\Ic(N))\to
\Hom_{\Hb'_\reg}(\Ic(M)_\reg,\Ic(N)_\reg)$$
is surjective.
An element $f\in\Hom_{\Hb'_\reg}(\Ic(M)_\reg,\Ic(N)_\reg)$
is a horizontal $W$-invariant section 
of the bundle $\Hom_{\Rb_\reg}(\Ic(M)_\reg,\Ic(N)_\reg)$ over $T_\reg$.
Given $\b\in\Delta_+$, we expand $f=\sum_{k\geq k_0}(1-z^\b)^kf_k$
locally near a generic point of $\{z^\b=1\}$,
with $f_k$ holomorphic on the divisor and $f_{k_0}$ not identically zero.
The residue of the connection on $\{z^\b=1\}$ is constant and has eigenvalues
$0,\pm 2h_{0\b}$, see 3.4. 
Thus $k_0\geq 0$ since $2h_{0\b}\notin\ZZ$.
\qed

\vskip3mm

\noindent{\bf Remark.}
Observe that $\Pu(\mu_0)^\nabla\neq\Pu(\e^{\mu_0})$ in general.
For instance, in type $A_1$, if $\l_0=\rho/2$ and $h_0=1/2$
then $\e^{s_\heartsuit\l_0}=\ell_0^{-1}$, and
$\Pu(s_\heartsuit\l_0)^\nabla=\Pu(\ell_0)\neq\Pu(\ell_0^{-1})$.
See 6.2 for more details.

\subhead 4.2\endsubhead
We do not know how to compute $\Pu(\mu_0)^\nabla$ for all
$\mu_0\in\hat W\l_0$. However we can prove a parabolic analogue
to Theorem 4.1.$(i)$ which is sufficient to recover the category
$\Oc'$ in type $A$, see Section 5.

Fix a non-empty subset $J\subsetneq\hat I$.
The group $W_J$ acts on $\Hb'_J$ on the right by translations.
The quotient is a left $\Hb'_J$-module 
which is naturally identified with $\Sb'$.
Let $O'\subset\hat W\l_0$ be a finite subset
such that $W_JO'=O'$.
The proof of Lemma 2.1 gives
$u[O']\subseteq\sum_{u'\leq u}[O'] u'$ 
for all $u\in W_J$.
Hence the ideal $[O']\subset\Sb'$ is preserved by $\Hb'_J$.
Set $P_J(O')=\Hb'\otimes_{\Hb'_J}\Sb_{O'}$,
and $1_{O'}=1\otimes 1\in P_J(O')$. 
The module $P_J(O')$
lies in ${}^{\l_0}\Oc$,
and is generated by $1_{O'}$ over $\Hb'$ 
with the defining relations 
$[O']\, 1_{O'}=0$ and $W_J\, 1_{O'}=1_{O'}.$
If $J\subseteq I$ then $P_J(O')=\Ic(\Pu_J(O'))$,
where $\Pu_J(O')=\Hbu'\otimes_{\Hb'_J}\Sb_{O'}$.

From now on we assume that $J\subseteq I$.
Set $C_{J,+}=\{\mu_0\in X_\RR\,;\,
(\mu_0\,:\,\a_j^\vee)=0,\,
(\mu_0\,:\,\a_k^\vee)>0,\,\forall j\in J,\,k\notin J\}$.
There is a unique representation of $\Hb_J$ on $\Sb$ such that
$t_j\, 1=\zeta_{0j}$ and $\Sb$ acts by multiplication. 
Set $\underline{[E]}=\bigcap_{m\in E}\underline{[m]}$ 
and $\Sb_E=\Sb/\underline{[E]}$
for any subset $E\subseteq W\ell_0$.
If $O\subset W\ell_0$ is a $W_J$-orbit,
the ideal $\underline{[O]}\subset\Sb$ is preserved by $\Hb_J$.
Set $\Pu_J(O)=\Hbu\otimes_{\Hb_J}\Sb_O$ and $1_O=1\otimes 1$.
The $\Hbu$-module $\Pu_J(O)$ lies in ${}^{\ell_0}\Oc$,
and is generated by $1_O$ over $\Hbu$
with the defining relations 
$\underline{[O]}\, 1_O=0$ and $t_j\, 1_O=\zeta_{0j}1_O$
for each $j\in J.$

\proclaim{Proposition}
(i)
$P_J(O')$ is projective in ${}^{\l_0}\Oc$.

(ii)
If $m_0\in W\ell_0$, $\mu_0\in\hat W\l_0$ are such that $\e^{\mu_0}=m_0$
and $\mu_0\in W\l_0-\kappa$ with $\kappa\in Y$ far enough inside $C_{J,+}$,
then $\Pu(W_J\mu_0)^\nabla=\Pu(W_Jm_0)$.
\endproclaim

\noindent{\sl Proof.}
Set $M_{O'}=\{m\in M\,;\,[O']m=0\}$ for each $\Hb'$-module $M$.
By Lemma 2.1 the subspace $M_{O'}\subset M$ is preserved by $W_J$. 
Moreover $\Hom_{\Hb'}(P_J(O'), M)=(M_{O'})^{W_J}$.
Hence $P_J(O')$ is projective in ${}^{\l_0}\Oc$,
because the functor $M\mapsto\bigl(M_{O'}\bigr)^{W_J}$ 
from ${}^{\l_0}\Oc$ to vector spaces is exact.
Claim $(i)$ is proved.

The proof of $(ii)$ is the same as the proof of Theorem 4.1.$(i)$,
to which we refer for notations and details.
Set $\mu_0\in\hat W\l_0$, $O'=W_J\cdot\mu_0,$ and $O=W_J\cdot m_0$. 

\vskip2mm

\noindent $(A)$
We first prove that $\Pu_J(O')^\nabla_\CC$ is cyclic over $\Hbu_\CC$.
To do so, we deform $\Pu_J(O')_\CC$.
From now on $\k=\A,\F$ or $\CC$,
$w,w'\in W$, $v,v'\in W^J$, $u\in W_J$, 
and $\nu_0,\nu'_0\in O'$.

Put $\mu=\mu_0+\varpi\epsilon$, with $\epsilon\in X_0$.
Set $Q=W_J\hat W_{\mu_0}\cdot\mu$ and 
$\bar\nu_0=Q\cap(\nu_0+\varpi X_\CC)$.
Let $\Sb_{Q,\A}$, $\Sb_{\bar\nu_0,\A}$ be as in 2.1. 
Set $\Sb_{Q,\k}=\k\otimes_\A\Sb_{Q,\A}$ and  
$\Sb_{\bar\nu_0,\k}=\k\otimes_\A\Sb_{\bar\nu_0,\A}$ if $\k=\F,\CC$.
The ring $\Sb_{\bar\nu_0,\A}$ is local.
Let $\mb_{\bar\nu_0,\A}$ be the maximal ideal. 
We have $\Sb_{\bar\nu_0,\CC}=\Sb_{\nu_0,\CC}$, see 4.1$(A)$.
Thus $\Sb_{Q,\CC}=\Sb_{O',\CC}$,
because $\Sb_{Q,\A}=\Oplus_{\nu_0}\Sb_{\bar\nu_0,\A}$.
Let $\nu$ denote any element in $Q$.
The embedding $\Sb_{Q,\A}\subseteq\Oplus_\nu\Sb_{\nu,\A}$ is 
generically invertible. 

The $\Hb'_J$-action on $\Sb'$
descends to $\Sb_Q$ because $W_JQ=Q$.
Set $\Pu=\Hbu'\otimes_{\Hb'_J}\Sb_Q$.
The module $\Pu$ lies in $\Ocu'$ and $\Pu_\CC=\Pu_J(O')_\CC$.
Set $\psi_v=\phi'_v\otimes 1\in\Pu$.
Assume that $W_{\mu_0}\subseteq W_J$.
We claim that $\Pu=\Oplus_v\psi_v\cdot\Sb_Q$.
It is enough to prove it for $\k=\A$.
Recall that 
$\psi_v\in v\psi_1\cdot\pi_v+\Sum_{v'<v}v'\psi_1\cdot\Sb_{Q,\A}$.
The image of $\pi_v$ in $\Sb_{Q,\CC}$ is invertible,
because $\Sb_{Q,\CC}=\bigoplus_{\nu_0}\Sb_{\nu_0,\CC}$
and $\pi_v\notin\mb_{\nu_0,\CC}$
(indeed, 
$\Delta_{J,+}\subset\Delta_+^\vee\cap v^{-1}\!\Delta_+^\vee$ 
because $v\in W^J$, hence
$\xi_{\a^\vee}(\nu_0)\neq 0$ if 
$\a^\vee\in\Delta_+^\vee\cap v^{-1}\!\Delta_-^\vee$ 
since $W_{\nu_0}\subseteq W_J$).
Therefore $\pi_v$ is invertible in $\Sb_{Q,\A}$.
The claim follows.

If $k\in\ZZ$ is non-zero the element
$t=k+{}^{v^{-1}}\!\!\xi_j-{}^{{v'}^{-1}}\!\!\xi_j$ 
is invertible in $\Sb_{Q,\F}$,
because $\Sb_{Q,\F}=\Oplus_\nu\Sb_{\nu,\F}$ 
and the projection of $t$ in $\Sb_{\nu,\F}$ 
is invertible (since $\epsilon\in X_0$).
Thus there is a unique fundamental matrix solution 
$G\,:\,\CC^I_\reg\to\End(\Pu_\F)$ 
of the trigonometric Knizhnik-Zamolodchikov connection of the form $G=Hz^{A_0}$,
with $H$ holomorphic on $\CC^I\setminus D_\Delta$ and $H(0)=\Id$. 
It yields a $\F$-linear isomorphism $\Pu^\nabla_\F\to\Pu_\F$.
From now on we identify the $\F$-vector spaces $\Pu_\F$, $\Pu_\F^\nabla$.
The $B_{\hat W}$-action on $\Pu^\nabla$ factorizes through
$\Hbu$ by Theorem 4.1.$(ii)$.
Thus $\Pu_\F$ admits left actions of $\Hbu'$ and $\Hbu$,
such that $y_j=\e^{\xi_j}$.
Moreover $\Pu^\nabla_\A\subset\Pu_\F$ is a $\Hbu_\A$-submodule,
and the canonical map $\CC\otimes_\A\Pu^\nabla_\A\to\Pu_\CC^\nabla$ 
is an isomorphism of $\Hbu_\CC$-modules.

We now fix $\mu_0$ as in $(ii)$.
Hence
$$(\nu_0\,:\,\b^\vee)\in\RR_{\ll 0}+i\RR,\qquad
\forall\b^\vee\in\Delta_+^\vee\setminus\Delta_{J,+}^\vee,\ \forall\nu_0.
\leqno(4.2.1)$$
In particular, $W_{\mu_0}\subseteq W_J$.
Assume that $s_jv>v$ and $s_jv\notin vW_J$. Hence $s_jv\in W^J$. 
We claim that
$$\forall p\in\Sb_{Q,\A},\ \exists x\in\Hbu_\A
\ \roman{such\ that}\ x\,\psi_v\in\psi_{s_jv}\cdot p
+\Sum_{v'<s_jv}\psi_{v'}\cdot\Sb_{Q,\F}.
\leqno(4.2.2)$$
We have $\Sb_\A\psi_v=\psi_v\cdot\Sb_{Q,\A}$,
because $y_j$ acts as $\e^{\xi_j}$ in $\Pu_\F$ 
and the ring homomorphism $\Sb_\A\to\Sb_{Q,\A}$, $p\mapsto p(\e^\xi)$
is surjective.
Therefore it is sufficient to prove that
$$t_j\psi_v\in\psi_{s_jv}\cdot\Sb_{Q,\A}^\times+
\Sum_{v'<s_jv}\psi_{v'}\cdot\Sb_{Q,\F}.\leqno(4.2.3)$$
The same argument as for (4.1.5) implies that
$t_j\psi_v\in\Sum_{v'\leq s_jv}\psi_{v'}\cdot\Sb_{Q,\F}$.
Set $P_j=\Hb'_j\otimes_{\Sb'}\Sb_{vQ}$.
The ring $\Hb_j$ acts on $P_{j,\F}$ by monodromy as in 4.1$(D)$.
Fix $\varphi_v,\varphi_{s_jv}\in P_j$ as in $4.1(D)$.
Let $\theta_j\,:\,P_j\to\Pu$ be the $\Hb'_j$-linear embedding 
such that $\varphi_w\mapsto\psi_w$ if $w=v$ or $s_jv$.
Using $\theta_j$ as in $4.1(D)$ we are reduced to prove that
$$t_j\varphi_v\in\varphi_{s_jv}\cdot\Sb_{vQ,\A}^\times+
\varphi_v\cdot\Sb_{vQ,\F}.$$
An element in $\Sb_{vQ,\A}$ is invertible if and only if
its image in $\Sb_{v\nu,\A}$ is invertible for each $\nu$. 
The projection $\Sb_{vQ}\to\Sb_{v\nu}$ 
yields a $\Hb'_j$-linear map $P_j\to P_j(v\nu)$. 
Using this map we are reduced to prove (4.1.9) again.
We have $v^{-1}\a_j^\vee\in\Delta_+^\vee\setminus\Delta^\vee_{J,+},$ 
because $s_jv>v$ and $s_jv\in W^J$.
Hence $(v\nu_0\,:\,\a_j^\vee)\notin\{0,\pm h_0\}+\ZZ_{\geq 0}$ by (4.2.1).
The claim (4.2.2) follows.

We now prove that (4.2.2) implies that $\Pu_\CC^\nabla=\Hbu_\CC\psi_1$.
If $\kappa$ is far enough inside $C_{J,+}$ there is an open convex cone
$\Cc\subset X_\CC\setminus\{0\}$ (i.e. $x+y,tx\in\Cc$ for each
$x,y\in\Cc$ and $t\in\RR_{>0}$) containing $Y_+\setminus\{0\}$
such that $v\nu_0-v'\nu'_0\in\Cc$ 
for each $v>v'$ and each $\nu_0,\nu'_0$.
Given $\eta_0,\eta'_0\in X_\CC$ we write $\eta_0\succ\eta'_0$ if 
$\eta_0-\eta'_0\in\Cc$.
Then $v\nu_0\succ v'\nu'_0$ if $v>v'$ or $v\nu_0>v'\nu'_0$.

Fix a $\A$-basis $(s_{\nu_0,t})$ of $\Sb_{\bar\nu_0,\A}$ for each $\nu_0$.
Write $\psi_{v,\nu_0,t}$ for $\psi_vs_{\nu_0,t}$.
By Proposition 7.1 there is a $\A$-basis 
$(b_{v,\nu_0,t})$ of $\Pu_\A^\nabla$ such that
$$b_{v,\nu_0,t}\in\psi_{v,\nu_0,t}+
\sum_{v'\nu'_0<v\nu_0}\sum_{t'}\psi_{v',\nu'_0,t'}\cdot\F.\leqno(4.2.4)$$

We first prove that $\psi_1\in\Pu_\A^\nabla$.
Since $\psi_1$ is a $\A$-linear combination of the elements
$\psi_{1,\nu_0,t}$, it suffices to check that
$b_{1,\nu_0,t}=\psi_{1,\nu_0,t}$ for each $\nu_0,t$.
By (4.2.4) it suffices to check that 
$v'\nu_0'\not <\nu_0$ for each $\nu'_0$, $v'$.
If $v'\neq 1$ then $v'\nu_0'-\nu_0\in\Cc$, hence
$v'\nu_0'\not <\nu_0$ because $(-\Cc)\cap Y_+=\emptyset$.
If $v'=1$ then $\nu_0'-\nu_0\notin Y\setminus\{0\}$
because a direct computation, 
using $\hat W_{\l_0}\subseteq W$ and $\nu_0\in W\l_0-\kappa$,
yields $\hat W_{\nu_0}\cap(Y\rtimes W_J)\subseteq W_J$.

Given $\nu_0,t$ there is $x\in\Hbu_\A$ such that 
$x\,\psi_1\in \psi_{v,\nu_0,t}+\Sum_{v'<v}\psi_{v'}\cdot\Sb_{Q,\F}$ by (4.2.2) 
and an obvious induction on $\ell(v)$.
Then 
$$x\,\psi_1\in 
b_{v,\nu_0,t}+\sum_{v'\nu'_0\prec v\nu_0}\sum_{t'}b_{v',\nu'_0,t'}\cdot\A$$
because $x\psi_1\in\Pu_\A^\nabla$.
Therefore $\Pu_\A^\nabla=\Hbu_\A\psi_1$,
hence $\Pu_\CC^\nabla=\Hbu_\CC\psi_1$.

\vskip2mm

\noindent $(B)$
Next, we prove that there is a unique surjective $\Hbu_\CC$-linear map
$\Pu_J(O)_\CC\to\Pu_\CC^\nabla$ such that $1_{O}\mapsto\psi_1$.
To do so we must prove that $\underline{[O]}\psi_1=0$, 
and $t_j\psi_1=\zeta_{0j}\psi_1$ for all $j\in J$. 

We have $\hat W_{\nu_0}=x_\kappa^{-1}W_{\e^{\nu_0}}x_\kappa$ for each $\nu_0$,
because $\hat W_{\l_0}=W_{\ell_0}$ and $\nu_0\in W\l_0-\kappa$.
Thus the map $p\mapsto p(\e^\xi)$ yields a ring isomorphism
$\Sb_{\e^{\nu_0},\CC}\to\Sb_{\nu_0,\CC}$, see $4.1.(B)$.
We have also a bijection of $W_J$-sets $O\simeq O'$, because 
$W_J\cap W_{m_0}=W_J\cap W_{\mu_0}$
(since
$W_{m_0}=x_\kappa\hat W_{\mu_0}x_\kappa^{-1}$,
and
$W_J\cap x_\kappa\hat W_{\mu_0}x_\kappa^{-1}=
W_J\cap W_{\mu_0}$ because $W_J$ centralizes $x_\kappa$).
Hence the map $p\mapsto p(\e^\xi)$ yields a ring isomorphism 
$\Sb_{O,\CC}\to\Sb_{O',\CC}$.
Hence $\underline{[O]}\psi_1=0$, because $[O']\psi_1=0$.

Assume that $j\in J$, $v=1$.
Then $t_j\psi_1=\psi_1\cdot p$ with $p\in\Sb_{Q,\F}$.
We claim that $p=\zeta_{0j}$.
For each $\nu$ the subspace $\psi_1\cdot\Sb_{\{\nu,s_j\nu\}}\subset\Pu$ 
is preserved by $\Hb'_j$.
Thus we are reduced to a computation in
$\Sb_{\{\nu,s_j\nu\}}$ over $\F$.
The result follows from \cite{C1}.

There is a unique surjective $\Hbu_\CC$-linear map
$\Pu_J(O)_\CC\to\Pu_\CC^\nabla$ 
such that $1_{O}\mapsto\psi_1$.
It is invertible because both modules have the same dimension
(since $\Sb_{O,\CC}\simeq\Sb_{O',\CC}$).
\qed

\subhead 4.3\endsubhead
Fix an integer $n>0$, and a subset $J\subseteq I$.
Given finite subsets $O\subset W\ell_0$ and $O'\subset\hat W\l_0$ 
which are preserved by $W_J$, we put $\Sb_{O',n}=\Sb'/[O']^n$ and 
$\Sb_{O,n}=\Sb/\underline{[O]}^n$.
Set $\Pu_J(O')_n=\Hbu'\otimes_{\Hb'_J}\Sb_{O',n}$,
$P_J(O')_n=\Ic(\Pu_J(O')_n),$
and $\Pu_J(O)_n=\Hbu\otimes_{\Hb_J}\Sb_{O,n}$.
Clearly $P_J(O')_n\in{}^{\{\l_0\}}\Oc$ and
$\Pu_J(O)_n\in{}^{\{\ell_0\}}\Oc.$
For each integer $n>0$, let
${}^{\l_0}\Oc_n\subset{}^{\{\l_0\}}\Oc$
be the full subcategory consisting of the modules $M$ such that, 
for each $m\in M$, there is a finite subset $E\subset\hat W\l_0$
with $[E]^nm=0$.
For a future use we need the following extension of 4.1-2.

\proclaim{Proposition}
(i)
$P_J(O')_n$ is projective in ${}^{\l_0}\Oc_n$.

(ii)
If $m_0\in W\ell_0$, $\mu_0\in\hat W\l_0$ are such that $\e^{\mu_0}=m_0$
and $\mu_0\in W\l_0-\kappa$ with $\kappa\in Y$ far enough inside $C_{J,+}$,
then $\Pu_J(W_J\mu_0)_n^\nabla=\Pu_J(W_Jm_0)_n$.

(iii)
The map
$\Mc\,:\,\Hom_{\Hb'}(P_{J_1}(O'_1)_n,P_{J_2}(O'_2)_n)\to
\Hom_{\Hbu}(\Pu_{J_1}(O_1)_n^\nabla,\Pu_{J_2}(O_2)_n^\nabla)$
is bijective. 
\endproclaim

\noindent{\sl Proof.}
For each $\Hb'$-module $M$ we set $M_{O',n}=\{m\in M\,;\,[O']^nm=0\}$.
The functor $M\mapsto\{m\in M\,;\,[O']^nm=0\}^{W_J}$
is exact on ${}^{\l_0}\Oc_n$ and is represented by $P_J(O')_n$.
Thus $P_J(O')_n$ is projective in ${}^{\l_0}\Oc_n$.
Claim $(ii)$ is proved as in 4.2, replacing everywhere
$\Sb_Q$ by $\Sb_{Q,n}=\Sb'/[Q]^n$.
The map in $(iii)$ is injective by Lemma 3.4.$(i)$
because $P_{J_2}(O'_2)_n$ is free over $\Rb$.
Any projective and indecomposable module $N$ 
in ${}^{\l_0}\Oc_n$ is a direct summand
of a module $P_\emptyset(\mu_0)_n$ with 
$\mu_0\in\hat W\l_0$, see Proposition 2.2.$(ii)$.
Since $P_\emptyset(\mu_0)_n\in\Ic(\Ocu')$,
the functor $\Mc$ is fully faithful on the projective modules
in ${}^{\l_0}\Oc_n$.
Thus the map in $(iii)$ is also surjective.
\qed

\head 5. Type $A$ case\endhead 
\subhead 5.1\endsubhead
Let $G^\vee$ be the simple simply connected and connected linear group
whose weight lattice is $X$ and whose root system is $\Delta$.
Thus $T^\vee$ is a maximal torus in $G^\vee$.
Let $\gen^\vee$ be the Lie algebra of $G^\vee$ over $\CC$.

Given $h_0\in\CC$, $\l_0\in X_\CC$ we set
$\ell'_0=\e^{\l_0}$, $\zeta'_0=\e^{h_0}$, and
$${}^{\ell'_0}\Nc_{\zeta'_0}=\{x\in\gen^\vee\,;\,
x\ \roman{is\ nilpotent\ and\ } 
\ad(\ell'_0)(x)=\zeta'_0x\}.$$
Let ${}^{\ell'_0}H\subseteq G^\vee(\CC)$ be the centralizer of $\ell'_0$. 
The group ${}^{\ell'_0}H$ acts on ${}^{\ell'_0}\Nc_{\zeta'_0}$ by conjugation.

Given $u_0\in\CC^\times$ satisfying (2.5.1) 
we set $\zeta_0=\e^{u_0h_0}$, 
$\tau_0=\e^{u_0}$, 
$\ell_0=\e^{u_0\l_0}$, and
$${}^{\ell_0}\Nc_{\zeta_0,\tau_0}=\{x(\varpi)\in\gen^\vee\otimes_\CC\F\,;\,
x(\varpi)\ \roman{is\ nilpotent\ and\ } 
\ad(\ell_0)(x(\tau_0\varpi))=\zeta_0x(\varpi)\}.$$
Let $\tilde G^\vee(\F)$ be the Kac-Moody central extension of $G^\vee(\F)$,
and
$${}^{\ell_0}H_{\tau_0}=\{g(\varpi)\in\tilde G^\vee(\F) \,;\,
\ad(\ell_0)(g(\tau_0\varpi))=g(\varpi)\}.$$
Note that ${}^{\ell_0}\Nc_{\zeta_0,\tau_0}\subseteq\gen^\vee\otimes_\CC\R$ 
and ${}^{\ell_0}H_{\tau_0}\subseteq\tilde G^\vee(\R)$,
where $\R=\CC[\varpi,\varpi^{-1}]$,
because $\tau_0$ is not a root of unity. 
The group ${}^{\ell_0}H_{\tau_0}\rtimes\CC^\times$ 
acts on ${}^{\ell_0}\Nc_{\zeta_0,\tau_0}$ :
the first factor acts by conjugation, the second by `rotation of the loops'.

\proclaim{Lemma}
The map $ev\,:\,\gen^\vee\otimes_\CC\R\to\gen^\vee$, $x(\varpi)\mapsto x(1)$
factorizes through a bijection 
${}^{\ell_0}\Nc_{\zeta_0,\tau_0}/({}^{\ell_0}H_{\tau_0}\rtimes\CC^\times)\to
{}^{\ell'_0}\Nc_{\zeta'_0}/{}^{\ell'_0}H$.
\endproclaim

\noindent{\sl Proof:}
We first claim that $ev$ restricts to an isomorphism 
${}^{\ell_0}\Nc_{\zeta_0,\tau_0}\to{}^{\ell'_0}\Nc_{\zeta'_0}$.
Given $x(\varpi)\in\gen^\vee\otimes_\CC\R$ we fix a decomposition  
$x(\varpi)=\sum_ix_i\otimes\varpi^{k_i}$, with $x_i\in\gen^\vee$, 
such that $x_i$ has the weight $\b_i^\vee$
and the elements $x_i\otimes\varpi^{k_i}$ 
are linearly independent over $\CC$. 
Using (2.5.1) we get
$$\ad(\ell_0)(x(\tau_0\varpi))=\zeta_0x(\varpi)
\iff (\l_0\,:\,\b_i^\vee)+k_i=h_0,\ \forall i.$$
In particular
$$\ad(\ell_0)(x(\tau_0\varpi))=\zeta_0x(\varpi)
\Rightarrow
\ad(\ell'_0)(x(1))=\zeta'_0x(1).$$
On the other hand,
if $\ad(\ell'_0)(x)=\zeta'_0x$
and $x=\sum_ix_i$ with $x_i$ of weight $\b_i^\vee$ and
$\b^\vee_i\neq\b^\vee_j$ if $i\neq j$, then
for each $i$ there is an integer $k_i$
such that $(\l_0\,:\,\b_i^\vee)+k_i=h_0$.
Thus the element 
$x(\varpi)=\sum_ix_i\otimes\varpi^{k_i}$ 
satisfies
$$\ad(\ell_0)(x(\tau_0\varpi))=\zeta_0x(\varpi)
\quad\roman{and}\quad x(1)=x.$$

If $\ad(\ell_0)(x(\tau_0\varpi))=\zeta_0x(\varpi)$,
$\ad(\ell_0)(y(\tau_0\varpi))=\zeta_0y(\varpi)$
and $x(1)=y(1)$, then, given decompositions
$x(\varpi)=\sum_ix_i\otimes\varpi^{k_i}$,
$y(\varpi)=\sum_jy_j\otimes\varpi^{\ell_j}$ 
as above, we get $\sum_ix_i=\sum_jy_j$,
and $k_i=\ell_j$ whenever the weights of $x_i$, $y_j$ coincide.
Thus $x(\varpi)=y(\varpi)$.

Obviously $x(1)$ is nilpotent if the element $x(\varpi)$ is nilpotent.
Conversely assume that $x(1)\in{}^{\ell'_0}\Nc_{\zeta'_0}$
and $\ad(\ell_0)(x(\tau_0\varpi))=\zeta_0x(\varpi)$.
Given $n>0$ we set 
$y(\varpi)=\ad(x(\varpi))^n\in\End(\gen^\vee\otimes_\CC\R)$.
Fix a decomposition 
$y(\varpi)=\sum_iy_i\otimes\varpi^{k_i}$,
such that $y_i\in\End(\gen^\vee)$ has the weight $\gamma_i^\vee$
and the elements $y_i\otimes\varpi^{k_i}$ are linearly independent over $\CC$. 
We have $(\l_0\,:\,\gamma_i^\vee)+k_i=nh_0$ for all $i$.
In particular $k_i=k_j$ whenever $\gamma_i^\vee=\gamma_j^\vee$.
Thus the operators $y_i$ are also linearly independent. 
Hence, if $y(1)=0$ then $y(\varpi)=0$.
Thus $x(\varpi)\in{}^{\ell_0}\Nc_{\zeta_0,\tau_0}$.
The claim is proved.

Given $x(\varpi)\in{}^{\ell_0}\Nc_{\zeta_0,\tau_0}$
the orbit $\ad({}^{\ell_0}H_{\tau_0})(x(\varpi))$ 
is a cone because $x(\varpi)$ is nilpotent  
(use the Jacobson-Morozov theorem as in 
Claim 2 in the proof of \cite{V2, Proposition 6.3} for instance).
Fix $k_i$, $\b^\vee_i$ as above. 
We have $(\l_0\,:\,\b_i^\vee)+k_i=h_0$ for all $i$. 
Hence for each $t\in\CC$ we get 
$x(\e^t\varpi)=\e^{h_0t}\ad(\e^{-t\l_0})(x(\varpi)).$
Clearly, $\e^{-t\l_0}\in {}^{\ell_0}H_{\tau_0}.$
Thus $x(z\varpi)\in\ad({}^{\ell_0}H_{\tau_0})(x(\varpi))$ 
for each $z\in\CC^\times$, i.e.
each ${}^{\ell_0}H_{\tau_0}$-orbit in 
${}^{\ell_0}\Nc_{\zeta_0,\tau_0}$
is preserved
by the action of $\CC^\times$ by rotation.
Therefore 
${}^{\ell_0}\Nc_{\zeta_0,\tau_0}/({}^{\ell_0}H_{\tau_0}\rtimes\CC^\times)=
{}^{\ell_0}\Nc_{\zeta_0,\tau_0}/{}^{\ell_0}H_{\tau_0}$.

Let $Z\subseteq\tilde G^\vee(\F)$ 
be the kernel of the obvious projection $\tilde G^\vee(\F)\to G^\vee(\F)$.
Thus $Z\simeq\CC^\times$.
Obviously, we have $Z\subseteq {}^{\ell_0}H_{\tau_0}$.
The map $ev\,:\,G^\vee(\R)\to G^\vee(\CC)$, 
$g(\varpi)\mapsto g(1)$ restricts to an isomorphism 
${}^{\ell_0}H_{\tau_0}/Z\to{}^{\ell'_0}H$.
Namely, both groups are connected by \cite{V2, Lemma 2.13},
$ev$ restricts to an injection 
${}^{\ell_0}H_{\tau_0}/Z\to{}^{\ell'_0}H$ 
(see \cite{BEG, Proposition 5.13} for instance),
and $ev$ yields an isomorphism of the Lie algebras of 
${}^{\ell_0}H_{\tau_0}/Z$ and ${}^{\ell'_0}H$.
Therefore $ev$ yields a bijection 
${}^{\ell_0}\Nc_{\zeta_0,\tau_0}/{}^{\ell_0}H_{\tau_0}
\to{}^{\ell'_0}\Nc_{\zeta'_0}/{}^{\ell'_0}H$.
\qed

\subhead 5.2\endsubhead
Set $\k=\CC$.
Put $\Tbu=\Oplus_{J\subseteq I}\Hbu\otimes_{\Hb_J}\Sb$.
The quantized affine Schur algebra is the ring $\Scbu=\End_\Hbu(\Tbu)$.
The right $\Sb^W$-action on $\Tb$ commutes to the left $\Hbu$-action.
It yields a ring homomorphism $\Sb^W\to\Scbu$.
Given $\ell'_0\in T^\vee$,
let ${}^{\{\ell'_0\}}\Scu$ be the full subcategory of $\Scbu$-$mof$
consisting of the modules which are annihilated
by some power of $\la\ell'_0\ra^W$.
Note that $\Scbu$-$mof=\bigoplus_{\ell'_0}{}^{\{\ell'_0\}}\Scu$,
where $\ell'_0$ varies in a set of representatives 
of the $W$-orbits in $T^\vee$.

From now on let $\Delta$ be of type $A_{d-1}$.
Since the parameters $\zeta_{0i}$ (resp. $h_{0i}$) are all equal 
we omit the subscript $i$.
To keep track of the parameters, 
we will index the categories considered so far by $\zeta_0,\tau_0$, etc.

\proclaim{Lemma}
The number of simple objects in ${}^{\{\ell'_0\}}\Scu_{\zeta_0'}$
is not less than the number of ${}^{\ell'_0}H$-orbits in 
${}^{\ell'_0}\Nc_{\zeta'_0}$.
\endproclaim

\noindent{\sl Proof.}
We only sketch the proof because the arguments are standard.
For any quasi-projective $\SL_d(\CC)\times\CC^\times$-variety $X$,
let $\Kb(X)$ be the complexified Grothendieck group of 
$\SL_d(\CC)\times\CC^\times$-equivariant coherent sheaves on $X$.
Set $\Bb=\Kb(point)$.
Recall that $\Spec(\Bb)=(T^\vee/W)\times\CC^\times.$
Let ${}^{\ell'_0}\Kb(X)_{\zeta_0'}$ be the specialization 
of the $\Bb$-module $\Kb(X)$
at $(W\ell'_0,\zeta_0')\in(T^\vee/W)\times\CC^\times.$

Let $X$ (resp. $\dot X$) be the variety of $d$-steps flags in $\CC^d$
(resp. the complete flag variety in $\CC^d$).
Let $Z\subseteq T^*X\times T^*X$ 
(resp. $\dot Z\subseteq T^*\dot X\times T^*\dot X$,
$Y\subseteq T^*\dot X\times T^*X$)
be the corresponding Steinberg varieties.
We endow $\Kb(Z)$, $\Kb(\dot Z)$ with an associative unital 
$\Bb$-linear product as in \cite{CG}.
We endow $\Kb(Y)$ with the $(\Kb(\dot Z), \Kb(Z))$-bimodule structure as in
\cite{GRV}.
We have a ring isomorphism $\Hbu\simeq\Kb(\dot Z)$, 
and $\Tbu\simeq\Kb(Y)$ over $\Hbu$.
It gives rise to a ring homomorphism $\Kb(Z)\to\Scbu$
which is invertible generically over $\Spec(\Bb)$.
The specialization map from the Grothendieck group of
$\Kb(Z)\otimes_{\CC[\zeta^{\pm 1}]}\CC(\zeta)$-$mof$
to the Grothendieck group of $\Kb(Z)_{\zeta_0'}$-$mof$ is surjective. 
Thus the pull-back map from the Grothendieck group of
${}^{\ell_0'}\Scu_{\zeta_0'}$
to the Grothendieck group of ${}^{\ell'_0}\Kb(Z)_{\zeta_0'}$-$mof$
is also surjective.
The simple ${}^{\ell'_0}\Kb(Z)_{\zeta_0'}$-modules
are labelled by
${}^{\ell'_0}\Nc_{\zeta'_0}/{}^{\ell'_0}H$ following \cite{V1} 
(see the remark below Theorem 4
if $\zeta'_0$ is a root of unity).
We are done.
\qed

\subhead 5.3\endsubhead
Let $\zeta'_0$, $\zeta_0$, $\tau_0$ be as in 5.1. 

\proclaim{Theorem}
Assume that $h_0\in\CC\setminus(1/2)\ZZ$.

(i)
$\Oc'_{h_0}$ and $\Scu_{\zeta'_0}$
are equivalent.

(ii)
If $\ell_0$, $\ell'_0$ are as in 5.1 then
${}^{\{\ell_0\}}\Oc_{\zeta_0,\tau_0}$ and ${}^{\{\ell'_0\}}\Scu_{\zeta'_0}$
are equivalent.
\endproclaim

\noindent{\sl Proof.}
We first prove $(i)$.
Fix $\l_0\in X_\CC$.
We can assume that $\hat W_{\l_0}\subseteq W$. 
Namely
for any $\pi\in\Omega$ the pull-back by the automorphism $\pi$ of $\Hb'$
yields an equivalence of categories 
${}^{\{\pi\l_0\}}\Oc_{h_0}\to{}^{\{\l_0\}}\Oc_{h_0}$.
Thus ${}^{\{w\l_0\}}\Oc_{h_0}$ and 
${}^{\{\l_0\}}\Oc_{h_0}$ are equivalent whenever $w\in \tilde W$.
Observe also that, in type $A$, 
there is always an element $w\in\tilde W$ such that 
$\hat W_{w\l_0}\subseteq W$ by Lemma 1.3.
In the other hand, for each $\pi$ the categories
${}^{\{\pi\ell'_0\}}\Scu_{\zeta'_0}$ and ${}^{\{\ell'_0\}}\Scu_{\zeta'_0}$ 
are equivalent,
because $\pi\ell'_0=w_\pi(\ell'_0\cdot\e^{w_\pi^{-1}\omega_\pi})$,
the element $\e^{w_\pi^{-1}\omega_\pi}\in G^\vee$ is central,
and for any $\ell\in T^\vee$ and $\omega$ in the center of $G^\vee$ we have
${}^\ell\Hbu\simeq{}^{\ell\omega}\Hbu$ 
(use \cite{CG, Proposition 8.1.5} for instance).

We have $\hat W_{\ell_0}=\hat W_{\l_0}$ because 
$$\matrix
x_\b w(\ell_0)=\ell_0&\iff
(\tau_0\otimes\beta)(\tau_0\otimes w\l_0)=\tau_0\otimes\l_0\hfill\cr
&\iff\b+w\l_0=\l_0\hfill\cr
&\iff x_\b w(\l_0)=\l_0.\hfill
\endmatrix
$$
Thus $\hat W_{\ell_0}=W_{\ell'_0}$, hence it is generated by reflections.
Moreover, 
$$\a^\vee\in\Delta^\vee_{(\ell_0)}\iff
\e^{u_0(\l_0\,:\,\a^\vee)}\in\Gamma\iff
(\l_0\,:\,\a^\vee)\in\ZZ+\ZZ h_0$$
by (2.5.1).
Therefore, Proposition 2.5.$(iii)$ yields an equivalence
${}^{\{\ell_0\}}\Oc_{\zeta_0,\tau_0}\simeq{}^{\{\l_0\}}\Oc_{h_0}.$

Since there is an involution of $\Hb'$ taking $h_0$ to $-h_0$,
we may assume that $h_0\notin\QQ_{>0}$. 
Thus the pair $(\tau_0,\zeta_0)$ is regular according to the
terminology in \cite{V2, Definition 2.14}.
Hence, by \cite{V2, Theorem 7.6 and Lemma 8.1} the simple objects in
${}^{\{\ell_0\}}\Oc_{\zeta_0,\tau_0}$, hence in
${}^{\{\l_0\}}\Oc_{h_0}$, are labelled by 
${}^{\ell_0}\Nc_{\zeta_0,\tau_0}/({}^{\ell_0}H_{\tau_0}\rtimes\CC^\times)$.
In the following we construct a quotient functor 
${}^{\{\l_0\}}\Oc_{h_0}\to{}^{\{\ell'_0\}}\Scu_{\zeta'_0}$.
It is an equivalence by Lemmas 5.1 and 5.2.

For each integer $n>0$ we set
${}^{\ell'_0}\Tbu_n=\Tbu\otimes_\Sb\Sb/[\ell'_0]_W^n$ 
and ${}^{\ell'_0}\Scbu_n=\Scbu\otimes_\Sb\Sb/[\ell'_0]_W^n$.
Thus $\End_{\Hbu}({}^{\ell'_0}\Tbu_n)={}^{\ell'_0}\Scbu_n.$
Note that $[\ell'_0]_W=\underline{[W\ell'_0]}$
by the Pittie-Steinberg theorem,
because $W_{\ell'_0}$ is generated by reflections.
If $J\subseteq I$ then $\Sb_{W\ell'_0,n}=\Oplus_O\Sb_{O,n}$,
where $O$ is any $W_J$-orbit in $W\ell'_0$.
Hence ${}^{\ell'_0}\Tbu_n=\Oplus_{J\subseteq I}\Oplus_O\Pu_J(O)_n$.

According to Proposition 4.3,
for each $J,O$ we can fix a $W_J$-orbit $O'\subset\hat W\l_0$
such that $\Pu_J(O')_n^\nabla=P_J(O)_n$.
Set ${}^{\l_0}\Tb_n=\bigoplus_{O'}P_J(O')_n$.
Consider the functor
$\Mc\,:\,\Oc'_{h_0}\to\Ocu_{\zeta'_0}$
introduced in Section 3.
We have $\Mc({}^{\l_0}\Tb_n)={}^{\ell'_0}\Tbu_n$, 
${}^{\l_0}\Tb_n$ is projective in ${}^{\l_0}\Oc_{n,h_0}$,
and $\End_{\Hb'}({}^{\l_0}\Tb_n)={}^{\ell'_0}\Scbu_n$
by Proposition 4.3.$(i),(ii)$.
Thus we have the quotient functor
$$F_n\,:\,{}^{\l_0}\Oc_{n,h_0}\to{}^{\ell'_0}\Scbu_n-mof,\quad
M\mapsto\Hom_{\Hb'}({}^{\l_0}\Tb_n,M).$$
It is an equivalence because both 
categories have the same (finite) number of simple objects.

On the other hand 
$${}^{\{\l_0\}}\Oc_{h_0}=\ind_n{}^{\l_0}\Oc_{n,h_0},\quad
{}^{\{\ell'_0\}}\Scu_{\zeta'_0}=\ind_n
({}^{\ell'_0}\Scbu_n-mof),$$
where $\ind$ stands for the inductive limit of categories.
The functors $F_n$ are compatible with the
inductive systems of categories.
Consider the $\Hb'$-module ${}^{\l_0}\Tb_\infty=\pro_n{}^{\l_0}\Tb_n$.
Note that ${}^{\l_0}\Tb_\infty\notin\Oc'_{h_0}$, 
because the $\Sb'$-action is not locally finite.
The natural map
$F_n(M)\to\Hom_{\Hb'}({}^{\l_0}\Tb_\infty,M)$
is an isomorphism for each $M\in{}^{\l_0}\Oc_{n,h_0}$.
Hence the functor
$$F_\infty\,:\,{}^{\{\l_0\}}\Oc_{h_0}\to{}^{\{\ell'_0\}}\Scu_{\zeta'_0},\quad
M\mapsto\Hom_{\Hb'}({}^{\l_0}\Tb_\infty,M),$$ 
is an equivalence of categories.
Claim $(i)$ follows, because 
$$\Oc'_{h_0}=\bigoplus_{\l_0}{}^{\{\l_0\}}\Oc_{h_0},\quad
\Scu_{\zeta'_0}=\bigoplus_{\ell'_0}{}^{\{\ell'_0\}}\Scu_{\zeta'_0},
$$
where $\l_0$ (resp. $\ell'_0$)
runs in a set of representatives of $\hat W$-orbits in $X_\CC$
(resp. of $W$-orbits in $X_\CC/Y$) and
the map $\l_0\mapsto\ell'_0$ is a bijection
$X_\CC/\hat W\to(X_\CC/Y)/W$.

The proof of $(ii)$ follows immediately from the proof of $(i)$.
\qed

\head 6. Another example\endhead 
\subhead 6.1\endsubhead
For any $\Hb'$-module $M$ in ${}^{\l_0}\Oc$, the character of $M$ is
the element 
$$\ch(M)=\sum_{\mu\in\hat W\l_0}\dim(M_\mu)\,\eps^\mu\in\ZZ X_\CC,$$
where $M_\mu$ is as in 2.2.
We do not assume that the root system is of type $A$ anymore,
but we restrict our attention to one single block in $\Oc'$.
Let $n$ be the Coxeter number.
Fix a positive integer $k$ prime to $n$.
Put $h_{0i}=h_0=k/n\in\QQ$, $\l_0=\rho/n\in X_\QQ$,
$\zeta_{0i}=\zeta_0=\e^{k/n}$, and $\ell_0=\e^{\rho/n}$.
Note that $\hat W_{\l_0}=\{1\}$.
For any $j\in\ZZ$ we set 
$\Delta^\vee(j)=\{\b^\vee\in\Delta^\vee\,;\,(\rho\,:\,\b^ \vee)=j\}$.
Set $k=an+b$, with $0<b<n$.
We have
$$\Hc_{\l_0}=\{(\b^\vee,-a),(\g^\vee,-1-a)\,;\,
\b^\vee\in\Delta^\vee(-b),\,\g^\vee\in\Delta^\vee(n-b)\}.$$
For each non-empty subset 
$J\subseteq I_k:=\Delta^\vee(-b)\cup\Delta^\vee(n-b)$ we set 
$$A_J=\{\mu\in X_\RR\,;\,
(\mu\,:\,\b^\vee),(\mu\,:\,\g^\vee)-1<a,\,
\forall\b^\vee\in J\cap\Delta^\vee(-b),\,
\forall\g^\vee\in J\cap\Delta^\vee(n-b)\}.$$
The function $J\mapsto A_J$ is decreasing.
Put $D_J=A_J\setminus\bigcup_{J'\supsetneq J}\bar A_{J'}.$
The sets $D_J$ are the affine domains.

\proclaim{Lemma}
The simple objects $\{V_J\}$ in ${}^{\l_0}\Oc$ are uniquely labelled
by non-empty subsets $J\subseteq I_k$ in such a way that 
$$\ch(V_J)=\sum_{A_w\subseteq D_J}\eps^{w\l_0}.$$
\endproclaim

\noindent{\sl Proof.}
Fix $v_0\in\CC^\times$ not a root of unity, and set
$\zeta_0=(v_0)^k$, $\tau_0=(v_0)^n$, $\ell_0=v_0^{n\l_0}$.
By Proposition 2.5 the categories
${}^{\{\ell_0\}}\Oc$ and ${}^{\{\l_0\}}\Oc$ 
are equivalent. 
The simple modules in ${}^{\{\ell_0\}}\Oc$ 
are classified in \cite{V2}, 
and the Jordan-H\"older factors of induced modules are given there
via intersection cohomology of some stratified variety.
In our case, the corresponding variety is $\CC^{\hat I}$,
with the stratification induced by the coordinate hyperplanes.
This yields 
$$\sum_{J'\supseteq J}\ch(V_{J'})=\sum_{A_w\subseteq A_J}\eps^{w\l_0}.$$
\qed

\vskip3mm

For all $\mu_0\in\hat W\l_0$ we have 
$\ch(P(\mu_0))=\sum_{w\in\hat W}\eps^{w\l_0}$,
because $\hat W_{\l_0}=\{1\}$.
In particular $P(\mu_0)$ is indecomposable, because it is generated
by the one-dimensional subspace $P(\mu_0)_{\mu_0}$.
By the proposition above the modules $P(\mu_0)$
and $\bigoplus_J V_J$ are equal in the Grothendieck ring.
There are $2^{r+1}-1$ affine domains in $X_\RR$, 
where $r$ is the rank of $\gen^\vee$. 
The corresponding projective objects in ${}^{\l_0}\Oc$ 
are the projective covers $P_J$ of the simple modules $V_J$,
for each non-empty subset $J\subseteq I_k$.
The set $D_{I_k}$ is the unique bounded affine domain. 
We have $\Mc(V_{I_k})={0}$ because $V_{I_k}$ is finite-dimensional.

There are $2^{r+1}-2$ domains in $X_\RR$. 
The corresponding projective objects in ${}^{\ell_0}\Ocu$ are
the modules $\Mc(P_J)$ with $J\subsetneq I_k$ non-empty,
by Theorem 4.1.$(i)$.
We claim that $\Mc(P_{I_k})=\Pu_I(W\ell_0)$.
To prove the claim, observe that
$\Hom_{\Hb'}(P_I(W\l_0), V_{I_k})=(\Oplus_{\mu_0\in W\l_0}(V_{I_k})_{\mu_0})^W$.
Hence $P_I(W\l_0)$ surjects to $V_{I_k}$, because
$\Oplus_{\mu_0\in W\l_0}(V_{I_k})_{\mu_0}\neq\{0\}$ 
by the proposition above and $V_{I_k}$ is simple.
The module $P_I(W\l_0)$ is projective in ${}^{\l_0}\Oc$.
Hence it contains the projective cover of $V_{I_k}$ as a direct summand.
Thus $P_I(W\l_0)=P_{I_k}$, because $\ch(P_I(W\l_0))=\ch(P_{I_k})$.
On the other hand $\Pu_I(W\l_0)^\nabla=\Pu_I(W\ell_0)$ by Theorem 4.2
(with $J=I$). We are done.

Note that $\Mc(P_{I_k})=\Sb_{W\ell_0}$, and that
${}^{\ell_0}\Hbu=\oplus_{w\in W}\Pu(w\ell_0)$, hence
${}^{\ell_0}\Hbu$ is a sum (with positive multiplicities)
of the modules $\Mc(P_J)$ with $J\subsetneq I_k$.
Thus there is a quotient functor
${}^{\l_0}\Oc\to\End_{\Hbu}({}^{\ell_0}\Hbu\oplus\Sb_{W\ell_0})$-$mof$.
Therefore ${}^{\l_0}\Oc$ is equivalent to 
$\End_{\Hbu}({}^{\ell_0}\Hbu\oplus\Sb_{W\ell_0})$-$mof$,
because both categories have the same number of simple modules.
More generally, let ${}^{\{\ell_0\}}\Ccu$
be the full subcategory of
$\End_{\Hbu}(\Hbu\oplus\Sb)$-$mof$
consisting of the modules which are annihilated 
by some power of $\la\ell_0\ra^W$.

\proclaim{Proposition}
The category ${}^{\{\l_0\}}\Oc$ is equivalent to ${}^{\{\ell_0\}}\Ccu$.
\endproclaim

\subhead 6.2\endsubhead
We give more details in type $A_1$. 
Then $\l_0=\rho/2$, $h_0=1/2$, $\zeta_0=-1$, and $\ell_0=i\otimes\a_1$.
There are 3 simple objects $V(s_\heartsuit\l_0),
V(s_1\l_0),V(\l_0)$ in ${}^{\l_0}\Oc$, 
such that
$\ch(V(\l_0))=\eps^{\l_0},$ and
$$\ch(V(s_\heartsuit\l_0))=\sum_{j\in 1+4\ZZ_{<0}}(\eps^{j\l_0}+\eps^{-j\l_0}),\quad
\ch(V(s_1\l_0))=
\eps^{-\l_0}+\sum_{j\in 1+4\ZZ_{>0}}(\eps^{j\l_0}+\eps^{-j\l_0}).$$
The representation of $\Hb'$ on $V(\l_0)$ takes 
$\xi_1$ to $1/4$, and $s_1,s_\heartsuit$ to $1$.
The module $V(s_j\l_0)$ is the quotient of $\Hb'$ 
by the left ideal generated by $\{\xi_1-(s_j\l_0)_1, s_j+1\}$
for each $j=\heartsuit,1$. 
The modules $P(\l_0)$, $P(s_\heartsuit\l_0)$, $P(s_1\l_0)$
are the projective covers of 
$V(\l_0)$, $V(s_\heartsuit\l_0)$, $V(s_1\l_0)$
respectively in ${}^{\l_0}\Oc$.

There are 2 simple objects $\Vu(\ell_0)$, $\Vu(\ell_0^{-1})$ 
in ${}^{\ell_0}\Ocu$.
The module $\Vu(\ell_0^{\pm 1})$
is one-dimensional such that $t_1,y_1$ acts as $-1,\pm i$,
and $\Pu(\ell_0^{\pm 1})$ is the projective cover of 
$\Vu(\ell_0^{\pm 1})$ in ${}^{\ell_0}\Ocu$.

We have $\Mc(V(\l_0))=0$ because $\dim V(\l_0)<\infty$.
We have $\Mc(V(s_1\l_0))=\Vu(\ell_0^{-1})$ because
$V(s_1\l_0)$ is induced from the one-dimensional $\Hbu'$-module
such that $W$ acts via the signature, and
$\Vu(\ell_0^{-1})$ is the one-dimensional 
$\Hbu$-module such that $t_j$ acts by -1.
Similarly $\Mc(V(s_\heartsuit\l_0))=\Vu(\ell_0)$. 
Then, we get easily
$\Mc(P(s_\heartsuit\l_0))=\Pu(\ell_0)$ 
and $\Mc(P(s_1\l_0))=\Pu(\ell_0^{-1})$.
There is an exact sequence
$$0\to V(s_\heartsuit\l_0)\oplus V(s_1\l_0)\to P(\l_0)\to V(\l_0)\to 0.$$
It yields
$\Mc(P(\l_0))=\Vu(\ell_0)\oplus\Vu(\ell_0^{-1})$.
Set $O'=\{-\l_0,\l_0\}$ and $O=\{\ell^{-1},\ell_0\}$.
Then $\Vu(\ell_0)\oplus\Vu(\ell_0^{-1})=\Pu_I(O)$, and 
$P(\l_0)=P_I(O')$ because
the map $P(\l_0)\to P_I(O')$, $1_{\l_0}\mapsto(\xi_1+{1\over 4})1_{\pm\l_0}$
is surjective and $\ch P(\l_0)=\ch P_I(O')$.
Thus $\Mc(P_I(O'))=\Pu_I(O)$.
To conclude, note that
$\Sb_O=\Pu_I(O)$ and 
${}^{\ell_0}\Hbu=\Pu(\ell_0)\oplus\Pu(\ell_0^{-1}).$ 

\head 7. Appendix\endhead 
\subhead 7.1\endsubhead
Recall that $\A=\CC[[\varpi]]$, $\F=\CC((\varpi))$.
Fix a commutative $\A$-algebra $\Sb_\A$ which is free of rank $e$ over $\A$.
Let $(s_u)$ be a $\A$-basis of $\Sb_\A$.
Set $\Sb_\k=\k\otimes_\A\Sb_\A$ if $\k=\CC$ or $\F$.
Assume that $\Sb_\CC$ is a local Artinian ring with maximal ideal $\mb_\CC$.
Then $\Sb_\A$ is also a local ring.
Let $\mb_\A\subset\Sb_\A$ be the maximal ideal.
Let $V_\A$ be a free right $\Sb_\A$-module of rank $d$, with basis $(e_r)$.
From now on $r,s$ belong to $\{1,2,...d\}$,
and $u,v$ to $\{1,2,...e\}$.
We write $e_{ru}$ for $e_rs_u$.

Let $\nabla=\d-\sum_jA_j\d z_j/z_j$ 
be a linear integrable meromorphic connection
over $\CC^I$, with $A_j=\sum_{\b\geq 0}A_{j\b}z^\b$ 
and $A_{j\b}\in\End(V_\A)$.
The space of horizontal sections $V_\A^\nabla$ 
is a free $\A$-module of rank $de$.
Set $V_\k^\nabla=V_\A^\nabla\otimes_\A\k$.

Assume that 
$A_{j0}(e_{ru})=e_{ru}m_{rj}$ with $m_{rj}\in\Sb_\A$ such that 
$k+m_{rj}-m_{sj}\in\Sb_\F^\times$ for each integer $k\neq 0$.
Let $\mu_{rj}$ be the image of $m_{rj}$ in the residue field $\Sb_\A/\mb_\A$.
Set $m_r=\sum_j m_{rj}\otimes\a_j$ and
$\mu_r=\sum_j\mu_{rj}\otimes\a_j$. 

There is a unique fundamental matrix solution 
$G\,:\,\CC^I\setminus D_\infty\to\End(V_\F)$ of the form $G=Hz^{A_0}$, 
with $H\,:\,\CC^I\to\End(V_\F)$ holomorphic such that $H(0)=\Id$.
Set $f_{ru}=Ge_{ru}$. 
Then $(f_{ru})$ is a $\F$-basis of $V_\F^\nabla$.

There is an integer $k_0\leq 0$ such that
$f_{ru}\varpi^{-k_0}\in V_\A^\nabla$ for each $u,r$.
Put $\zeta_j=\log z_j$.
Let $V_\CC[\zeta]$ be the set of $V_\CC$-valued polynomials in the $\zeta_j$'s,
$W=V_\CC[\zeta][[\varpi]]\varpi^{k_0}$, and 
$W[[z]]$ be the set of $W$-valued formal series in the $z_j$'s.
Write $W[[z]]'\subset W[[z]]$ 
for the set of formal series without constant term.
Then $f_{ru}$ has an expansion in 
$e_{ru}z^{m_r}+W[[z]]' z^{\mu_r}.$

The following proposition is standard,
but we have not found a convenient reference.

\proclaim{Proposition}
There is a $\A$-basis $(b_{ru})$ in $V_\A^\nabla$
such that $b_{ru}\in f_{ru}+\sum_{\mu_s>\mu_r}\sum_v f_{sv}\F$.
\endproclaim

\noindent{\sl Proof.}
Note that $e_{ru}z^{m_r-\mu_r}\in W$ because $m_{rj}-\mu_{rj}\in\mb_\A$. 
Consider a formal series $b_{ru}=\sum_{\b\geq 0}b_{ru\b}z^{\mu_r+\b}$,
with $b_{ru\b}\in W$ and $b_{ru0}=e_{ru}z^{m_r-\mu_r}$.
It is the expansion of an horizontal section in $V_\A^\nabla$ if and only if
for all $j\in I$ we have
$$\partial_{\zeta_j}b_{ru\b}+b_{ru\b}(\b_j+\mu_{rj})-A_{j0}(b_{ru\b})
=\sum_{\g<\b}A_{j,\b-\g}(b_{ru\g}),
\ \forall\b\geq 0.
\leqno(7.1.1)$$
We have
$\partial_{\zeta_j}b_{ru0}+b_{ru0}\mu_{rj}-A_{j0}(b_{ru0})=0$ 
because $A_{j0}(b_{ru0})=b_{ru0}m_{rj}$.
Assume that $b_{ru\g}$ satisfies (7.1.1) for each $\g<\b$.
Recall that for all $j \in I$, $c\in W$ and $B\in\End(V_\A)$, 
there is an element $b\in W$ such that $\partial_{\zeta_j}b-B(b)=c$
(solve the equation term by term using asymptotic expansions of
$b,c,B$ in series in $\varpi$.
It is done inductively on the exponent of $\varpi$).
Hence, for each $j$ there is a non empty set of solutions 
$b_{ru\b}\in W$ to $(7.1.1)$.
There is a common solution for all $j$
because $\nabla$ is integrable.
Therefore, for each $(r,u)$ there is an
horizontal section $b_{ru}\in V_\A^\nabla$
with an expansion in $e_{ru}z^{m_r}+W[[z]]' z^{\mu_r}$.
These sections form a $\A$-basis of $V_\A^\nabla$ because
$(e_{ru})$ is a $\A$-basis of $V_\A$. 
Fix elements $x_{sv}\in\F$ such that 
$$b_{ru}-\sum_{s,v}f_{sv}x_{sv}=0.\leqno(7.1.2)$$
We must prove that $\mu_s>\mu_r$ if $x_{sv}\neq 0$ and $(s,v)\neq(r,u)$,
and that $x_{ru}=1$.

Consider expansions in $\varpi$ of the summands in (7.1.2).
Given $s$, let $\b_sz^{\mu_s}$ be the constant term in
$-\sum f_{sv}x_{sv}$ where the sum is over all $v$ such that
$(s,v)\neq(r,u)$,
and let $\a_rz^{\mu_r}$ the constant term in
$b_{ru}-f_{ru}x_{ru}.$
Then $\a_r,\b_s$ are holomorphic with asymptotic expansions
$\a_r(z),\b_s(z)$ in $V_\CC[\zeta][[z]]$.
Moreover the constant term $\b_s(0)\in V_\CC[\zeta]$ of the 
non zero series $\b_s$ are linearly independent.
Fix $\nu\geq 0$ minimal such that
$\a_rz^{\mu_r}=\g_rz^{\nu+\mu_r}$
and $\g_r$ has an asymptotic expansion in $V_\CC[\zeta][[z]]$ 
with non-zero constant term.
Then (7.1.2) gives
$$\g_rz^{\nu+\mu_r}+\sum_{s}\b_{s}z^{\mu_s}=0.\leqno(7.1.3)$$
We claim that $\nu>0$ and that there is an index $s$ such that 
$\b_s\neq 0$ and $\nu+\mu_r=\mu_s$.
Then, setting
$\g'_r=(\g_r+\sum_{\mu_s=\nu+\mu_r}\b_s)z^{\nu-\nu'}$ with $\nu'\geq\nu$
minimal such that $\g'_r(0)\neq 0$, and
$\b'_s=\b_s$ if $\mu_s\neq\nu+\mu_r$ and 0 else, (7.1.3) yields
$$\g'_rz^{\nu'+\mu_r}+\sum_{s}\b'_{s}z^{\mu_s}=0.$$
Once again there is an index $s$ such that 
$\b'_s\neq 0$ and $\nu'+\mu_r=\mu_s$.
By induction we have proved that
$\mu_s>\mu_r$ for each pair $(s,v)\neq(r,u)$ such that $x_{sv}\neq 0$.
Moreover $x_{ru}=1$ because $\nu>0$.
To prove the claim recall the following fact : 

\vskip2mm

\noindent (7.1.4)
given an equation $\Sum_{t=1}^mv_tz^{\nu_t}=0$
with $\nu_t\in X_\CC$ and $v_t$ holomorphic with an expansion 
$v_t(z)\in V_\CC[\zeta][[z]]$,
if the constant terms $v_t(0)$ are non-zero 
then $\nu_1,...\nu_m$ are not all different.

\vskip2mm

\noindent
(It is sufficient to prove this for $I=\{1\}$.
If $\nu_1,...\nu_m$ are all different we can fix
$\zeta\in\CC$ such that $|\e^\zeta|<1$ 
and $|\e^{\nu_1\zeta}|$,...$|\e^{\nu_m\zeta}|$ are distincts. 
Assume that $|\e^{\nu_{t_1}\zeta}|>|\e^{\nu_{t_2}\zeta}|>\cdots>
|\e^{\nu_{t_m}\zeta}|.$ 
Setting $\zeta\mapsto k\zeta$ with $k\gg 0$,
the equation $\Sum_{t=1}^mv_t(\e^{k\zeta})\e^{k\nu_t\zeta}=0$
yields $v_{t_1}(0)=0$). 
If $\nu=0$ then $x_{ru}\neq 1$.
Hence the elements $\gamma_r(0)$, $\b_s(0)$ with $s$ such that $\b_s\neq 0$
are linearly independent, and (7.1.3) yields a contradiction with (7.1.4).
The rest of the claim is immediate from (7.1.4) again.
\qed

\subhead 7.2\endsubhead
Let $\Ab$ be a ring with a unity,
and $S$ be an infinite (countable) set.
Put $\Ab^S=\bigoplus_{s\in S}\Ab$, 
and $\M_S(\Ab)=\Hom_\Ab(\Ab^S,\Ab^S)$
(with respect to the right $\Ab$-action on $\Ab$).
Elements in $\M_S(\Ab)$ may be viewed as infinite matrices 
whose columns have only finitely many entries.
If $\Ab$ is a topological ring we endow $\M_S(\Ab)$
with the finite topology :
a system of neighborhood of an element $f$ is formed by the subsets
$$\{f'\in\M_S(\Ab)\,;\,f(x)-f'(x)\in U^S,\,\forall x\in\Ab^E\},$$
where $E\subset S$ is finite and $U\subset\Ab$ is an open neighborhood of zero.
Recall that a $\Ab$-module $M$ is smooth if 
the annihilator in $\Ab$ of any element is open.
Let $\Ab$-$mod^\infty$ be the category of smooth finitely generated 
$\Ab$-modules.

\proclaim{Proposition}
The categories $\Ab$-$mod^\infty$ and $\M_S(\Ab)$-$mod^\infty$ 
are equivalent. 
\endproclaim

\noindent{\sl Proof :}
Set $\Bb=\M_S(\Ab)$.
The ring $B$ acts on $\Ab^S$ on the left, 
the ring $\Ab$ acts on $\Ab^S$ on the right.
To simplify assume that the topology on $\Ab$ is discrete.
The general case is identical.
We must prove that $\Ab$-$mod$ and $\Bb$-$mod^\infty$ 
are equivalent.
Consider the functors
$$\matrix
F\,:\,\Ab-\Mc od\to\Bb-\Mc od,&
M\mapsto\Ab^S\otimes_\Ab M,\hfill\cr\cr
G\,:\,\Bb-\Mc od\to\Ab-\Mc od,&
N\mapsto\Hom_{\Bb}(\Ab^S,N).\hfill
\endmatrix$$
The left $\Ab$-action on $G(N)$ comes from the right $\Ab$-action on $\Ab^S$.
The functor $G$ is exact because $\Ab^S$ is projective in $\Bb$-$\Mc od$.
The functor $F$ is obviously exact.

$(i)$
We have
$$GF(M)=\Hom_\Bb(\Ab^S,\Ab^S\otimes_\Ab M)=\Hom_\Bb(\Ab^S,\Ab^S)\otimes_\Ab M,$$
because $\Ab^S$ is finitely generated over $\Bb$.
The right $\Ab$-action on $\Hom_\Bb(\Ab^S,\Ab^S)$
comes from the right $\Ab$-action on $\Ab^S$.
Using commutation with elementary matrices, we get
$$\Hom_\Bb(\Ab^S,\Ab^S)=
\{f_a\,:\,\Ab^S\to\Ab^S,\,v\mapsto va\,;\,a\in\Ab\}\simeq
\Ab$$
(isomorphism of $(\Ab,\Ab)$-bimodules).
Thus $GF(M)=M.$

$(ii)$
The natural evaluation map
$$\phi_N\,:\,FG(N)=\Ab^S\otimes_\Ab\Hom_\Bb(\Ab^S,N)\to N$$
is a morphism of $\Bb$-modules.
We claim that $\phi_N$ is bijective if $N\in\Bb$-$mod^\infty$.

To prove the surjectivity it is sufficient to
assume that $N$ is smooth and cyclic.
For any finite set $E\subset S$, set 
$\Ib_E=\{f\in\Bb\,;\,f(x)=0,\,\forall x\in\Ab^E\}$.
Then it is enough to assume $N=\Bb/\Ib_E$,
because the ideals $\Ib_E$ form a basis of open neighborhoods of zero
in $\Bb$.
Clearly $\Bb/\Ib_E\simeq(\Ab^S)^E$ over $\Bb$.
Moreover $FG(\Ab^S)^E=F(\Ab)^E=(\Ab^S)^E$, by $(i)$,
and $\phi_N$ is the identity if $N=(\Ab^S)^E$.

We now prove the injectivity.
The exact sequence
$$0\to\Ker(\phi_N)\to FG(N)\to N\to 0$$
yields an exact sequence
$$0\to G(\Ker(\phi_N))\to G(N)\to G(N)\to 0$$
by (i), where the third map is 
$G(\phi_N)=\Id_{G(N)}$.
Thus $G(\Ker(\phi_N))=\{0\}.$
The $\Bb$-module $\Ker(\phi_N)$ is smooth, because $FG(N)$ is smooth.
Hence, for any finitely generated submodule $N'\subset\Ker(\phi_N)$
we have $G(N')=\{0\}$ and the map $\phi_{N'}$ is surjective.
Thus $N'=\{0\}.$
Therefore $\Ker(\phi_N)=\{0\}.$

$(iii)$
It is sufficient to check $G(\Bb$-$mod^\infty)\subset\Ab$-$mod$ 
on smooth cyclic $\Bb$-modules.
Thus it is enough to prove that
$G(\Bb/\Ib_E)\in\Ab$-$mod$ for each finite set $E\subset S$, see $(ii)$. 
This is obvious because
$G(\Bb/\Ib_E)\simeq G(\Ab^S)^E=\Ab^E$ by $(i)$.

$(iv)$
The inclusion $F(\Ab$-$mod)\subset\Bb$-$mod^\infty$ 
is obvious because $\Ab^S\subset\Bb$-$mod^\infty$.

\qed

\enddocument